%% file: main.tex
\documentclass{article}
\usepackage[utf8]{inputenc}

\usepackage[style=alphabetic,natbib,doi=false,isbn=false,url=false,eprint=false]{biblatex}
\usepackage{amsmath}
\usepackage{amssymb}
\usepackage{amsthm}
\usepackage{mathtools}

\usepackage{pgf,tikz,pgfplots}

\usepackage{tikz-cd}
\usetikzlibrary{patterns}

\usetikzlibrary{arrows,decorations.markings}

\usepackage{dirtytalk}

\usepackage{enumerate}
\usepackage[shortlabels]{enumitem}

\usepackage{subfigure}

\usepackage[margin=1in]{geometry}

\usepackage{xcolor}

\usepackage{setspace}
\theoremstyle{plain}
\newtheorem{thm}{Theorem}

\theoremstyle{plain}
\newtheorem{lemma}{Lemma}

\theoremstyle{plain}
\newtheorem{prop}{Proposition}

\theoremstyle{plain}
\newtheorem{cor}{Corollary}

\usepackage[T2A,T1]{fontenc}

\usepackage[symbol]{footmisc}

\usepackage{array}
\usepackage{makecell}

\renewcommand{\d}{\,\mathrm{d}}
\newcommand{\limn}{\lim_{n\rightarrow \infty}}

\newcommand{\tor}{\mathbb{T}^2}

\usepackage{subfiles}

\usepackage{multicol}
\providecommand{\keywords}[1]
{
  \small	
  \textbf{\textit{Keywords---}} #1
}

\providecommand{\Acknowledgements}[1]
{
  \small	
  \textbf{\textit{Acknowledgements---}} #1
}

\title{A Continuous Family of Non-Monotonic Toral Mixing Maps}

\author{J. Myers Hill$^{1\,\star}$, R. Sturman$^{2}$, M. C. T. Wilson$^{3}$  \\
        \small $^{1}$EPSRC CDT in Fluid Dynamics, University of Leeds, Leeds LS2 9JT, United Kingdom \\
        \small $^{2}$School of Mathematics, University of Leeds, Leeds LS2 9JT, United Kingdom \\
        \small $^{3}$School of Mechanical Engineering, University of Leeds, Leeds LS2 9JT, United Kingdom \\
        \small $^\star$ E: scjmh@leeds.ac.uk
}

\pgfplotsset{compat=1.17}

\begin{document}

\date{}

\maketitle

\begin{abstract}
 
We establish the mixing property for a family of Lebesgue measure preserving toral maps composed of two piecewise linear shears, the first of which is non-monotonic. The maps serve as a basic model for the `stretching and folding' action in laminar fluid mixing, in particular flows where boundary conditions give rise to non-monotonic flow profiles. The family can be viewed as the parameter space between two well known systems, Arnold's Cat Map and a map due to Cerbelli and Giona, both of which possess finite Markov partitions and straightforward to prove mixing properties. However, no such finite Markov partitions appear to exist for the present family, so establishing mixing properties requires a different approach. In particular we follow a scheme of Katok and Strelcyn, proving strong mixing properties with respect to the Lebesgue measure on two open parameter spaces. Finally we comment on the challenges in extending these mixing windows and the potential for using the same approach to prove mixing properties in similar systems.
 
%  \color{black}

\end{abstract}


\Acknowledgements{JMH is supported by EPSRC under Grant Ref. EP/L01615X/1.}

%\tableofcontents

\section{Introduction}
\subfile{sections/intro}

\section{Statement of Results}\label{sec:statement}
\subfile{sections/statement}

\section{Proof Outline}\label{sec:outline}

\subfile{sections/outline}

\section{Perturbation from Arnold's Cat Map}
\label{sec:eta}

\subsection{Establishing non-uniform hyperbolicity}

\label{sec:etaLyp}
\subfile{sections/etaPert/etaLyp}

\subsection{Establishing ergodicity}
\label{sec:etaMan}
\subfile{sections/etaPert/etaMan}

\subsection{Establishing the Bernoulli property}
\label{sec:etaRepMan}
\subfile{sections/etaPert/repeatedManifold}

\section{Perturbation from Cerbelli \& Giona's Map}
\label{sec:eps}

\subsection{Establishing non-uniform hyperbolicity}
\label{sec:epsLyp}
\subfile{sections/epsPert/epsLyp}

\subsection{Establishing ergodicity}
\subfile{sections/epsPert/epsMan}

\subsection{Establishing the Bernoulli property}
\subfile{sections/epsPert/epsRepMan}

\subfile{sections/thmProof}

\section{Final Remarks}\label{sec:final}
\subfile{sections/finalRemarks}

\section{Appendix}
\label{sec:appendix}
\subfile{sections/appendix/epsBoundsAnalysis}

\section{Supplementary Material}
\label{sec:suppMaterial}
\subsection{From section \ref{sec:eta}}

\subfile{sections/etaPert/supp}

\newpage

\subsection{From section \ref{sec:eps}}

\subfile{sections/epsPert/piecewiseLinearCurves}

\newpage

\printbibliography

% \bibliographystyle{alpha}
% \bibliography{refs.bib}

\end{document}

%% file: sections/intro.tex
Two-dimensional measure-preserving discrete-time dynamical systems are both rich in behaviour and relevant to a wide variety of applications. For example, as stroboscopic maps of fluid flow they constitute a model of kinematic mixing~\cite{ottino1989kinematics}; as canonical examples of Hamiltonian systems such as forced pendulums or kicked rotators \cite{ott2002chaos}; as fundamental models in fast dynamo theory \cite{childress1995stretch} and quantum chaos \cite{keating1991cat}. The richness of the dynamical behaviour can be seen in the observations that the dynamics may be integrable, but also may exhibit chaotic behaviour. That is, within two-dimensional maps, hyperbolicity is compatible with area-preservation, allowing access to the complete ergodic hierarchy, including ergodicity, measure-theoretic mixing, the Bernoulli property, etc.

This richness can be illustrated by considering the family of maps given by the transformation $H:(x,y) \to (x',y')$ of the 2-torus $\mathbb{T}^2$ into itself, given by
\begin{eqnarray}
x' &=& x+ f(y) \\
y' &=& y+x'.
\end{eqnarray}
Interpreting $H$ as the composition of a pair of shears $H = G \circ F$, with $F(x,y) = x+f(y)$, $G(x,y) = y+x$ clarifies that Lebesgue measure is preserved by $H$, regardless of the choice of $f$. In the case of the Cat Map, $f(y)=y$ imposes a constant, hyperbolic, Jacobian at every point in $\mathbb{T}^2$.  This fact provides the means to establish immediately dynamical properties, such as unstable manifolds all lying in the same direction, a positive Lyapunov exponent for every trajectory, and ergodic properties, such as strong mixing, the Bernoulli property and exponential decay of correlations. The uniform hyperbolicity of the Cat Map might be a desirable property, but is also strong enough to preclude many applications. 

The strict condition of uniformity of the hyperbolicity may be broken in a number of ways. A typical method is to slow down the expansion of tangent vectors. The first such example of a non-uniformly hyperbolic $C^{\infty}$ area-preserving map on $\mathbb{T}^2$ was the Katok map \cite{katok1979bernoulli}, in which trajectories near the hyperbolic fixed point at the origin are slowed down, with that fixed point becoming neutral. This is sufficient to produce zero Lyapunov exponents for some trajectories (although at almost every initial condition these remain non-zero), and thus non-hyperbolicity. In spite of the loss of uniform hyperbolicity, exponential decay of correlations are retained \cite{pesin2019thermodynamics}. 

Another example which breaks the uniformity of expansion is a linked twist map. Defined on a subset of $\mathbb{T}^2$ we replace function $f(y)$ with a piecewise smooth, non-decreasing function $\hat{f}(y)$, such that $d\hat{f}/dy =0$ over some sub-interval of $[0,1]$.  Now unstable leaves are oriented in a continuum of directions, but, crucially, all contained in the positive quadrant of tangent space, which makes the demonstration of the mixing property relatively straightforward.  Such a map retains the Bernoulli property of the Cat Map \cite{przytycki_ergodicity_1983}, but the rate of mixing is slowed to polynomial decay of correlations \cite{sturman2013rate,springham_polynomial_2014}.

One more example destroying the simplicity of the Cat Map can be found in the discontinuous sawtooth map. Here $f(y) = Cy$, with $C>0$, so that $C=1$ recovers the Cat Map. When $K$ is any other positive integer the map is continuous on the torus and the same analysis applied. When $K$ is non-integer however, the map becomes discontinuous, and although stable and unstable manifolds exist locally almost everywhere, these may be arbitrarily short, cut up by the dense countable set of discontinuity lines created by iterating the map. Nevertheless, the map retains its ergodicity~\cite{vaienti_ergodic_1992} as the parameter $C$ is perturbed from an integer.
		
For all the above examples, the map could be described as monotonic, in the sense that $f(y)$ is non-decreasing in each case. Much more complicated dynamics is possible if this condition is broken, as can be seen in the rich behaviour of the Chirikov--Taylor Standard Map \cite{chirikov1971research}. This well-known map, for which $f(y) = \frac{K}{2\pi} \sin 2 \pi y$, where $K$ is a parameter, can exhibit co-existence of invariant circles, elliptic islands and chaotic seas, due to the lack of an invariant cone in tangent space. The wider range of possible directions for unstable leaves allows for the possibility of expansion being immediately counteracted in the following iterate, and the consequent failure of hyperbolicity. 

A piecewise--linear version of the standard map was studied in \cite{wojtkowski1981model,bullett1986invariant}, where $f(y) = K\left( | y- 1/2| - 1/4 \right)$, and shown for certain parameter values to be non-uniformly hyperbolic ($K\geq 4$) and mixing ($K>K_0 \approx 4.0329$). For $K<4$ the map admits both chaotic and elliptic invariant domains; mixing properties over such a chaotic domain are shown for the $K=1$ map in \cite{liverani1995ergodicity}.  A different piecewise--linear adaptation of the standard map is that introduced by Cerbelli and Giona \cite{cerbelli_continuous_2005}, and proposed as a ``continuous archetype of area-preserving non-uniform chaos''. This map takes $f(y) = 2y$ if $y \in [0,1/2]$ and $f(y) = 2(1-y)$ if $y \in [1/2,1]$.  Like the Cat Map, the Cerbelli-Giona map (hereafter CG map) has a finite Markov partition~\cite{mackay_cerbelli_2006}, and so only a finite number of possible directions for piecewise linear segments in the unstable and stable leaves. 

Various generalisations to the CG map have been proposed, for example, in \cite{demers_family_2009} a family of maps designed to preserve the Markov structure is examined, while in \cite{mackay_cerbelli_2006} a number of perturbations preserving the pseudo-Anosov nature of the map are proposed. A smooth perturbation was considered in \cite{cerbelli_characterization_2008} and dynamical properties such as topological entropy were studied numerically, but the mixing property was not demonstrated. Here we take $f(y) = y/(1-\eta)$ if $y \in [0,1-\eta]$ and $f(y) = (1-y)/\eta$ if $y \in [1-\eta,1]$. At $\eta=0$ this gives the Cat Map, at $\eta=1/2$ the CG map, and at $\eta=1$ the map is periodic with period 6. We focus on the parameter space between the Cat map and the CG map, demonstrating the Bernoulli property over two subsets of $0<\eta<1/2$. In section \ref{sec:statement} we state our results and these subsets explicitly, while in section \ref{sec:outline} we summarise the steps the proof requires. Section \ref{sec:eta} deals with the parameter range near the Cat Map, and section \ref{sec:eps} with parameters near the CG Map. To keep the argument concise we move three derivations from section $\ref{sec:eps}$ to the appendix, section \ref{sec:appendix}. We conclude with some final remarks in section \ref{sec:final}.Explicit expressions for certain coordinates used in sections \ref{sec:eta}, \ref{sec:eps} are given as supplementary material, section \ref{sec:suppMaterial}.

%% file: sections/statement.tex
\begin{figure}[ht]
    \centering
    \begin{tikzpicture}
    \node at (-2.4,0) {
    \begin{tikzpicture}
    \draw (0,0) rectangle (4,4);
\draw (0,0)--(4,3)--(0,4);
\draw[<->] (4.1,4)--(4.1,3);
\node[scale=1.5] at (4.3,3.5) {$\eta$};
\draw[->] (0.4,2.8)--(2,2.8);
\node[scale=1.7] at (2,4.3) {$F$};
    \end{tikzpicture}
    };
    \node at (2.4,0) {
    \begin{tikzpicture}
    \draw (0,0) rectangle (4,4);
\draw (0,0)--(4,4);
\draw[->] (2.8,0.4)--(2.8,2);
\node[scale=1.7] at (2,4.3) {$G$};
    \end{tikzpicture}
    };
    \end{tikzpicture}
    \caption{A family of area preserving maps $H = G \circ F$ parameterised by $\eta>0$. Taking $\eta = 0$ gives the Cat map, taking $\eta=\frac{1}{2}$ gives \citeauthor{cerbelli_continuous_2005}'s map, both with well understood mixing properties.}
    \label{fig:mapDefn}
\end{figure}
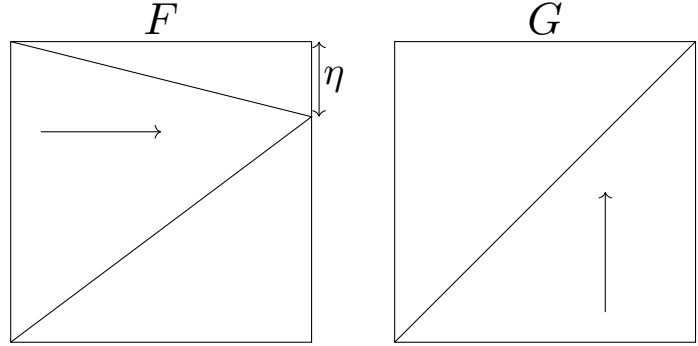

We consider the Lebesgue measure preserving map $H:\tor \rightarrow \tor$, taken as the composition of two orthogonal shears $H=G \circ F$, shown in Figure \ref{fig:mapDefn}. Taking local coordinates $(x,y) \in (\mathbb{R} / \mathbb{Z})^2$, $F$ maps
\[ (x,y) \mapsto 
\begin{cases}
\left(  x + \frac{1}{1-\eta} y , y  \right) \text{ mod 1 } & \text{ for } y \leq 1-\eta \\

\left(  x + \frac{1}{\eta} (1-y) , y  \right) \text{ mod 1 } & \text{ for } y \geq 1-\eta \\
\end{cases}\]
and $G$ maps $(x,y) \mapsto (x,y+x) \text{ mod 1}$, where $\eta$ is some real parameter $0 \leq \eta \leq \frac{1}{2}$. Note that $H$ is piecewise linear, with derivative
\[DH_1 = \begin{pmatrix} 1 & \frac{1}{1-\eta} \\1 & \frac{2-\eta}{1-\eta} \end{pmatrix} \]
for $0<y<1-\eta$, and
\[ DH_0 = \begin{pmatrix} 1 & -\frac{1}{\eta} \\1 & \frac{\eta -1}{\eta} \end{pmatrix}  \] for $1-\eta<y<1$. $DH$, then, is defined everywhere but the set $\mathcal{D} = \{(x,y) \,|\, y \in \{0,1-\eta \}\}$. The inverse map $H^{-1} = F^{-1} \circ G^{-1}$ is differentiable outside of the set $\mathcal{D}' = G(\mathcal{D})$. 

The aim of this paper is to prove mixing properties for $H$ over a wide parameter range. In particular, we prove:

\begin{thm}
\label{thm:Bernoulli}
$H$ has the Bernoulli property for $0<\eta<\eta_1$ and $\eta_2 \leq \eta<\eta_3$ where
$\eta_1 \approx 0.324$, $\eta_2 \approx 0.415$, and $\eta_3 \approx 0.491$.
\end{thm}

%% file: sections/outline.tex
Our scheme for proving the Bernoulli property is to satisfy the qualifications given in the following theorem from \cite{katok_invariant_1986}, paraphrased in \cite{sturman_ottino_wiggins_2006}.

\begin{thm}[\citeauthor{katok_invariant_1986}]
\label{thm:katok-strelcyn}
Let $f: X \rightarrow X$ be a measure preserving map on a measure space $(X,\mathcal{F},\mu)$ such that $f$ is $C^2$ smooth outside of a \emph{singularity set} $S$ where differentiability fails. Suppose that the Katok-Strelcyn conditions hold:
\begin{enumerate}[label={\bfseries (KS\arabic*):}]
    \item $\exists \, a,C_1>0$ s.t. $\forall \, \epsilon>0$, $\mu(B_\varepsilon(S)) \leq C_1 \varepsilon^a$.
    \item $\exists \, b,C_2>0$ s.t. $\forall \, z \in X \setminus S$, $||D^2_zf|| \leq C_2 \,d(z,S)^{-b}$ where $D^2_zf$ is the second derivative of $f$ at $z$.
    \item Lyapunov exponents exist and are non-zero almost everywhere.
\end{enumerate}
Then at almost every $z$ we can define local unstable and stable manifolds $\gamma_u(z)$ and $\gamma_s(z)$. Suppose that the manifold intersection property holds:
\begin{enumerate}[label={\bfseries (M):}]
    \item For almost any $z,z'\in X$, $\exists \, m,n$ s.t. $f^m(\gamma_u(z)) \cap f^{-n}(\gamma_s(z')) \neq \varnothing$.
\end{enumerate}
Then $f$ is ergodic. Furthermore the Bernoulli property holds, provided we can show the repeated manifold intersection property:
\begin{enumerate}[label={\bfseries (MR):}]
    \item For almost any $z,z'\in X$ we can find $M,N$ such that for all $m>M$ and $n>N$, $f^m(\gamma_u(z)) \cap f^{-n}(\gamma_s(z')) \neq \varnothing$.
\end{enumerate}
\end{thm}

The scheme extends Pesin theory (establishing ergodic properties of $C^2$ smooth non-uniformly hyperbolic systems, \citealp{pesin_characteristic_1977}) to systems which are smooth outside of some singularity set. The conditions \textbf{(KS1-2)} ensure that this set has manageable influence, and follow easily from our map's definition. We take our map as $f = H$, our domain as $X = \tor$, and our singularity set as $S = \mathcal{D}$. Taking $\mu$ to be the Lebesgue measure on $\tor$, clearly $\mu(S)=0$. When we say `for almost any $z \in \tor$', we will be referring to the full measure set $X' = \tor \setminus S_\infty$, $S_\infty = \bigcup_{k \geq 0} H^{-k}(\mathcal{D}) \cup \bigcup_{k \geq 0} H^{k}(\mathcal{D}')$, where $H$ and all its powers $H^k$, $k\in \mathbb{Z}$ are differentiable. Since we can cover $\mathcal{D}$ with arbitrarily thin rectangles, \textbf{(KS1)} follows for some $C_1>0$ with $a=1$. Since $H$ is piecewise linear, \textbf{(KS2)} follows trivially.

Moving onto \textbf{(KS3)}, we define the (forwards-time) Lyapunov exponent at a point $z \in \mathbb{T}^2$ in direction $v\in \mathbb{R}^2$ by
\[ \chi(z,v) = \limn \frac{1}{n} \log||DH^n_z v||, \]
where \[ DH^n_z = DH_{H^{n-1}(z)}\cdot ... \cdot DH_{H(z)} \cdot DH_z \]
is well defined at almost every $z$. We define $\log^+(\cdot) = \max \{ \log(\cdot), 0 \}$ and let $||\cdot||_{\mathrm{op}}$ be the operator norm. Existence of Lyapunov exponents almost everywhere follows from Oseledets' theorem (\citealp{oseledets_multiplicative_1968}) provided that $\log^+||DH||_{\mathrm{op}}$ is integrable. This clearly holds, so 
our first substantial task is proving that these exponents are non-zero. A particular form of Oseledets' theorem in two dimensions is useful here. We paraphrase from \cite{viana_lectures_2014}:

\begin{thm}[\citeauthor{oseledets_multiplicative_1968}, \citeauthor{viana_lectures_2014}]
\label{thm:Oseledets-2d}
Let $F:X \times \mathbb{R}^2 \rightarrow X \times \mathbb{R}^2$ be given by $F(x,v) = (f(x),A(x)v)$ for some measure preserving map $f$ on a 2-dimensional manifold $X$ and some measurable function $A: X \rightarrow \mathrm{GL}(2)$. Suppose $ \log^+|| A^{\pm 1} ||$ are integrable and define
\[ \lambda_+(x) = \limn \frac{1}{n} \log||A^n(x)||, \quad \lambda_-(x) = \limn \frac{1}{n} \log|| (A^n(x))^{-1}||^{-1}, \]
where $A^n(x) = A(f^{n-1}(x))\cdot ... \cdot A(f(x)) \cdot A(x)$. Then for almost every $x\in X$,
\begin{enumerate}
    \item either $\lambda_-(x) =\lambda_+(x)$ and
    \[ \limn \frac{1}{n} \log || A^n(x)v || = \lambda_\pm(x) \quad \forall v \in \mathbb{R}^2 \setminus \{ 0 \} \]
    \item or $\lambda_+(x) >\lambda_-(x)$ and there exists a vector line $E_x^s \subset \mathbb{R}^2$ such that
    \[ \limn \frac{1}{n} \log || A^n(x)v || = \begin{dcases}
    \lambda_-(x) & \text{for } v \in E_x^s \setminus \{0\}, \\
     \lambda_+(x) & \text{for } v \in \mathbb{R}^2 \setminus E_x^s.
     \end{dcases}
    \]
\end{enumerate}

\end{thm}

\begin{cor}
\label{cor:OurCorollary}
Further assuming that $A$ takes values in $\mathrm{SL}(2)$ gives $\lambda_-(x) = -\lambda_+(x)$. Hence if at some $x$ there exists $v_0 \in \mathbb{R}^2$ with $\lim_n \frac{1}{n} \log || A^n(x)v_0 || \neq 0$, it follows that $\lim_n \frac{1}{n} \log || A^n(x)v || \neq 0$ for all non-zero vectors $v$.
\end{cor}

Applying this corollary to the cocycle generated by the derivative $DH$ of our map $H$ gives an efficient scheme for establishing non-zero Lyapunov exponents. We let $A^n(z) = DH_z^n$, which takes values in SL(2). If there exists $v_0$ such that $|| DH^n_z v_0 ||$ grows exponentially with $n$, Corollary \ref{cor:OurCorollary} gives $\chi(z,v) \neq 0$ for all $v \neq 0$. Letting $\varepsilon = \frac{1}{2}-\eta$, we can either consider our system as an $\varepsilon$-perturbation from \citeauthor{cerbelli_continuous_2005}'s map, or as an $\eta$-perturbation from Arnold's Cat map. There is subset of the parameter space $\frac{1}{3}\leq\eta< \frac{1}{8}(9-\sqrt{33}) \approx 0.407$ in which island structures appear, splitting the parameter space into two sides. Proving \textbf{(M)} for the Cerbelli-Giona side follows a very similar argument to the Cat map side, but the calculations are generally more involved. For this reason we will begin by considering the $0<\eta<\frac{1}{3}$ perturbation, then continue with the $\varepsilon$-perturbation in section \ref{sec:eps}.

%% file: sections/etaPert/etaLyp.tex
In \cite{cerbelli_continuous_2005} a three element $ABC$ partition of the domain was defined with $H(A) \subset A \cup B$, $H(B) = C$, and $H(C)\subset A$. Their derivative matrix $DH|_A = DH_1$ was hyperbolic which, together with the fact that orbits leaving $A$ return after exactly two iterations in $A^c$, allowed Cerbelli and Giona to reduce much of the dynamics to that of a hyperbolic toral automorphism, with well understood mixing properties.

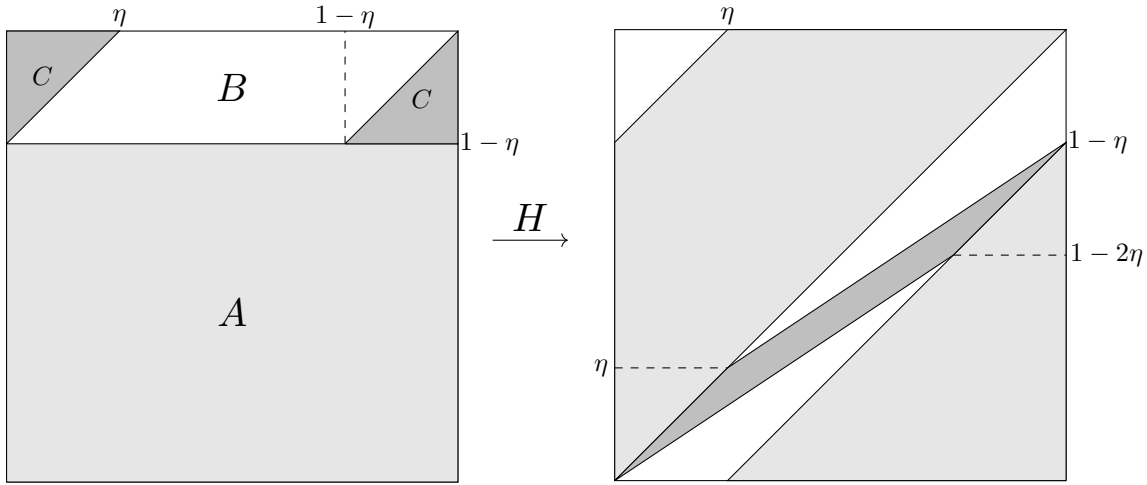
\begin{figure}[h]
    \centering
   \begin{tikzpicture}
\node at (-4,4) {
    \begin{tikzpicture}[scale = 0.6]
    \draw (0,0) rectangle (10,10);
    \draw[fill=gray!20] (0,0) -- (0,7.5) -- (10,7.5) -- (10,0) -- (0,0);
    \node[scale=1.5] at (5,3.75) {$A$};
    
    \draw[fill=gray!50] (7.5,7.5) -- (10,10) -- (10,7.5) -- (7.5,7.5);
    \node at (10.7,7.5) {$1-\eta$};
    \node[scale=1.5] at (5,8.75) {$B$};
    
    \draw[fill=gray!50] (0,7.5) -- (2.5,10) -- (0,10) -- (0,7.5);

    \node at (2.5,10.3) {$\eta$};
    \node at (7.5,10.3) {$1-\eta$};
    \draw[dashed] (7.5,7.5) -- (7.5,10);
    \node at (0.8,9) {$C$};
    \node at (9.2,8.5) {$C$};
     
    \end{tikzpicture}

};

\draw[->] (-1,4) -- (0,4);
\node[scale=1.5] at (-0.5,4.3) {$H$};

\node at (4,4) {
    \begin{tikzpicture}[scale = 0.6]
    \draw (0,0) rectangle (10,10);
    \draw[fill=gray!20] (0,0) -- (0,7.5) -- (2.5,10) -- (10,10) -- (0,0);
    \draw[fill=gray!20] (10,0) -- (10,7.5) -- (2.5,0) -- (10,0);
    
    \draw[fill=gray!50] (0,0) -- (2.5,2.5) -- (10,7.5) -- (7.5,5) -- (0,0);
    \node at (10.7,7.5) {$1-\eta$};
    \node at (10.9,5) {$1-2\eta$};
    \draw[dashed] (7.5,5) -- (10,5);
        \node at (2.5,10.3) {$\eta$};
   \draw[dashed] (0,2.5) -- (2.5,2.5);
   \node at (-0.3,2.5) {$\eta$};
    \end{tikzpicture}
};

\end{tikzpicture}    
    \caption{Partition of the torus for $H$, establishing return times to $A$ in $\{ 1,2,3 \}$. Case illustrated $\eta = \frac{1}{4}$, the image of the partition is also shown with consistent shading. \label{fig:returnTimePartitions}}
\end{figure}

While this approach is not possible for our family of maps, we do retain an upper bound on return times to $A$, illustrated by the partition of the domain given in 
Figure \ref{fig:returnTimePartitions}. One can show that $H(A) \subset A \cup B$, $H(B) \subset A \cup C$, $H(C) \subset A$ so that orbits leaving $A$ return after spending one or two iterations in $B\cup C$. We call the path an orbit takes around this partition its \emph{itinerary}. Any itinerary, for example
\[ AABCABAABABCA \dots \]
can be split up into itinerary blocks $I_j$ ending in $A$. In the above example this would look like
\[ A\quad A\quad BCA\quad BA\quad A \quad BA\quad BCA \dots.\]
There are three\footnote{Four if you include $CA$, the first block in the itinerary of a point starting in $C$, but this also has corresponding matrix $M_2$.} unique itinerary blocks
\[ I_1 = A, \quad I_2 = BA, \quad I_3 = BCA, \]
with corresponding matrices
\[ M_1 = DH_1, \quad M_2 = DH_1 \, DH_0, \quad M_3 = DH_1 \,DH_0^2. \]
Each $M_j$ is hyperbolic for $\eta$ strictly less than $\frac{1}{3}$, where $M_3$ loses hyperbolicity. Our parameter range, then, is $0<\eta<\frac{1}{3}$.

\begin{prop}
\label{prop:etaLyp}
We have non-zero Lyapunov exponents $\chi(z,v) \neq 0$ for almost every $z\in \tor $, $v\neq 0$, when $0<\eta<\frac{1}{3}$.
\end{prop}

\begin{proof}

Let $v$ be a non-zero vector in the tangent space at $x$. As the orbit starting at $x$ completes an itinerary block $I_j$, the effect on $v$ is to premultiply by the matrix $M_j$. Our aim is to find a vector $v_0$ which sees expansion in its norm after each itinerary block. The issue we have to overcome is the possibility that expansion by one matrix may be immediately undone by contraction from another. We do this by constructing an invariant, expanding cone.

We define a \emph{cone} $\mathcal{C}$ as a subset of $\mathbb{R}^2\setminus\{ 0\}$ such that if $v \in \mathcal{C}$ then $kv \in \mathcal{C}$ for any real $k \neq 0$.
Given a matrix $M$ we say that $\mathcal{C}$ is \emph{invariant} if $M\mathcal{C} \subset \mathcal{C}$. That is, vectors in the cone remain in the cone when premultiplied by $M$. We say that the cone is \emph{expanding} if $||Mv|| > ||v||$ for every $v\in \mathcal{C}$, where $||\cdot||$ is some norm we choose to put on the tangent space. In the tangent space take coordinates $(v_1,v_2)^T \in \mathbb{R}^2$. Since the transformations we are considering are linear and cones are double sided, the gradient of a vector is the only important feature.

Starting with invariance, if the gradients $g_j^u$, $g_j^s$ of the unstable, stable eigenvectors of $M_j$ satisfy
\[ g_1^s(\eta) < g_2^s(\eta) < g_3^s(\eta) < g_3^u(\eta) < g_2^u(\eta) < g_1^u(\eta),\]
then the cone bounded by (and including) the unstable eigenvectors of $M_1$ and $M_3$ will be invariant. Explicit expressions for these gradients will be given as supplementary material, and the chain of inequalities is easily verified for all $0<\eta<\frac{1}{3}$.

It is clear, then, that it is possible to construct an invariant cone and, in fact, we have multiple options. The minimal cone is the smallest gradient range we can take to include all the unstable eigenvectors, defined at each parameter value. This will be a particularly useful construction later on as it gives good bounds on the gradients of local unstable manifolds. Its $\eta$-dependence, however, makes the expansion factor calculations quite tedious. Given that $g_3^s(\eta) <  \inf_\eta g_3^u(\eta)$
across $0<\eta<\frac{1}{3}$, the cone bounded by (but not including) the vectors
$v_{\pm}$ with gradients $g^+ = \sup_\eta g_1^u(\eta) = \frac{2}{\sqrt{5}-1}$ and $g^- = \inf_\eta g_3^u(\eta) =1$ is invariant. Write this $\eta$-independent cone as $\overline{\mathcal{C}}$.

We will now show that $\overline{\mathcal{C}}$ is expanding. If across $0<\eta<\frac{1}{3}$ each of the $M_j$ expands both of the bounding vectors $v_\pm$, then the same holds for all vectors in the cone. To see this, note that (by hyperbolicity) $M_j$ expands its unstable eigenvector $v_u$, and contracts its stable eigenvector $v_s$. Let $\mathrm{ex}(v) := \frac{||M_jv||}{||v||}$, then $\mathrm{ex}(v_u)>1$ and $\mathrm{ex}(v_s)<1$. As we rotate $v$ from $v_u$ to $v_s$, we pass through one of $v_\pm$ and $\mathrm{ex}(v)$ has at most one local minimum. If $\mathrm{ex}(v_\pm)>1$, then this minimum must lie between $v_\pm$ and $v_s$, i.e. outside of the cone, so $\{ \mathrm{ex}(v) \,|\, v \in \overline{\mathcal{C}} \}$ is minimal at one of its boundaries. To simplify the calculations take $||\cdot||$ to be the $||\cdot||_\infty$ norm then $||(v_1,v_2)^T|| = |v_2|$ for all vectors in the cone, since within $\overline{\mathcal{C}}$ we always have $|v_2| \geq |v_1|$. Normalise the cone boundaries as $v_\pm = \left(\frac{1}{g^\pm},1 \right)^T$, now we can calculate:
\begin{multicols}{2}
\begin{itemize}

    \item $ ||M_1(1,1)^T|| = \frac{2 \eta - 3}{\eta - 1} > 3 $
    \item $ ||M_2(1,1)^T|| = \frac{3 \eta^{2} - 7 \eta + 3}{\eta \left(1-\eta\right)} > \frac{9}{2}$
    \item $ ||M_3(1,1)^T|| = 4 - \frac{10}{\eta} + \frac{3}{\eta^{2}} > 1$
    \item $ ||M_1\left(\frac{\sqrt{5}-1}{2}, 1 \right)^T||  = \frac{\left(1 + \sqrt{5}\right) \left(\eta - 1\right) - 2}{2 \left(\eta - 1\right)} > \frac{3+ \sqrt{5}}{2}$
    \item $ ||M_2\left(\frac{\sqrt{5}-1}{2},1 \right)^T||  = \frac{2 \sqrt{5} \eta^{2} - 3 \sqrt{5} \eta - 5 \eta + 6}{2 \eta \left(1 - \eta\right)} > \frac{39-7\sqrt{5}}{4} $
    \item $ ||M_3\left(\frac{\sqrt{5}-1}{2},1 \right)^T|| = \frac{ \left(3 \sqrt{5} -1 \right) \eta^{3} - \left(7 \sqrt{5} +7 \right) \eta^{2} + \left(3 \sqrt{5} + 17\right) \eta - 6}{2 \eta^{2} \left(\eta - 1\right)} > \frac{31-9\sqrt{5}}{4} $

\end{itemize}
\end{multicols}
for all $0 < \eta < \frac{1}{3}$, so that the cone is expanding across the parameter range.

\end{proof}

This establishes $H$ as non-uniformly hyperbolic over $0 < \eta < \frac{1}{3}$. The aim of the next section is to show that \textbf{(M)} holds, establishing ergodicity.

%% file: sections/etaPert/etaMan.tex
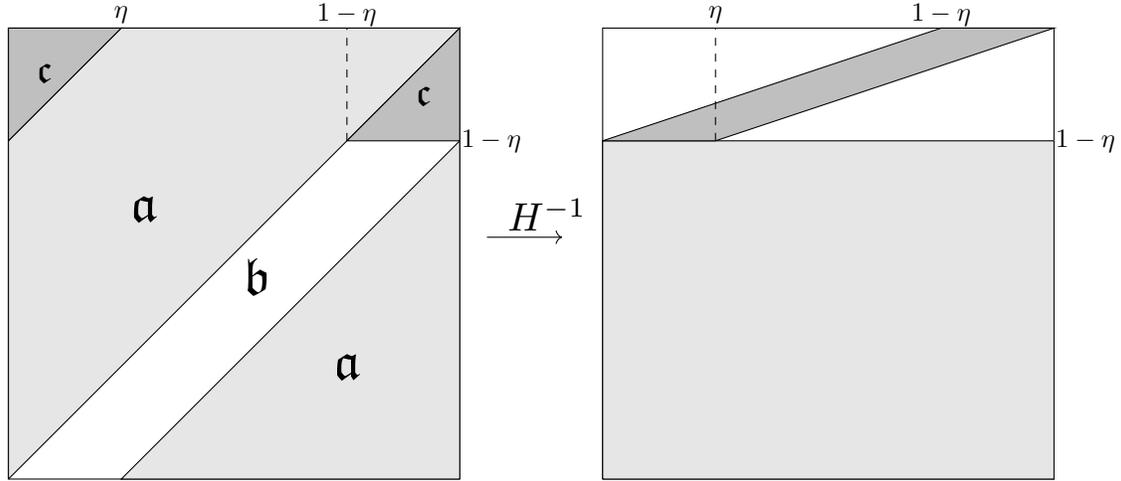
\begin{figure}
    \centering
      \begin{tikzpicture}

\node at (-3.9,-4) {

    \begin{tikzpicture}[scale = 0.6]
    \draw (0,0) rectangle (10,10);
    \draw[fill=gray!20] (0,0) -- (0,7.5) -- (2.5,10) -- (10,10) -- (0,0);
    \draw[fill=gray!20] (10,0) -- (10,7.5) -- (2.5,0) -- (10,0);
    
    \draw[fill=gray!50] (7.5,7.5) -- (10,10) -- (10,7.5) -- (7.5,7.5);
    \node at (10.7,7.5) {$1-\eta$};
    \draw[fill=gray!50] (0,7.5) -- (2.5,10) -- (0,10) -- (0,7.5);
    
    \node at (7.5,10.3) {$1-\eta$};
    \draw[dashed] (7.5,7.5) -- (7.5,10);
    \node[scale=1.5] at (0.8,9) {$\mathfrak{c}$};
    \node[scale=1.5] at (9.2,8.5) {$\mathfrak{c}$};
    \node[scale=2] at (3,6) {$\mathfrak{a}$};
    \node[scale=2] at (7.5,2.5) {$\mathfrak{a}$};
    \node[scale=2] at (5.5,4.5) {$\mathfrak{b}$};
    \node at (2.5,10.3) {$\eta$};
    
    \end{tikzpicture}
};

\draw[->] (-1,-4) -- (0,-4);
\node[scale=1.5] at (-0.2,-3.7) {$H^{-1}$};

\node at (4,-4) {
    \begin{tikzpicture}[scale = 0.6]
    \draw (0,0) rectangle (10,10);
    \draw[fill=gray!20] (0,0) -- (0,7.5) -- (10,7.5) -- (10,0) -- (0,0);
    
    \draw[fill=gray!50] (0,7.5) -- (2.5,7.5) -- (10,10) -- (7.5,10) -- (0,7.5);;
    \node at (10.7,7.5) {$1-\eta$};
    \node at (7.5,10.3) {$1-\eta$};
    \draw[dashed] (2.5,7.5) -- (2.5,10);
    \node at (2.5,10.3) {$\eta$};

\end{tikzpicture}
};

    \end{tikzpicture}
    \caption{A partition of the torus based on returns to $\mathfrak{a}$ under $H^{-1}$ and its image under $H^{-1}$. Case illustrated $\eta = \frac{1}{4}$.}
    \label{fig:backwardsPartition}
\end{figure}

Return time partitions and invariant cones can be similarly constructed for $H^{-1}$. These are useful for the next step, so we will give them now. Figure \ref{fig:backwardsPartition} shows the partition for returns to the set $\mathfrak{a}$. The itinerary blocks are follow the same pattern: $\mathfrak{a}$, $\mathfrak{b}\mathfrak{a}$, and $\mathfrak{b}\mathfrak{c}\mathfrak{a}$, with corresponding matrices $\mathfrak{M}_1$, $\mathfrak{M}_2$, and $\mathfrak{M}_3$ respectively. The eigenvectors of each of these matrices allow us to construct an invariant expanding cone $\mathcal{C}'$. Let $\mathfrak{g}_j^s(\eta)$, $\mathfrak{g}_j^u(\eta)$ be the gradients of the stable, unstable eigenvectors of $\mathfrak{M}_j$. One can verify that
\[ \mathfrak{g}_1^u(\eta) <  \mathfrak{g}_2^u(\eta) <  \mathfrak{g}_3^u(\eta) <  \mathfrak{g}_3^s(\eta) <  \mathfrak{g}_2^s(\eta) <  \mathfrak{g}_1^s(\eta)  \]
for $0< \eta < \frac{1}{3}$ so that we can take our minimal backwards cone to be the cone bounded by (and containing) the unstable eigenvectors of $\mathfrak{M}_1$ and $\mathfrak{M}_3$. As before, taking the union of these cones over $0<\eta<\frac{1}{3}$ gives an $\eta$-independent invariant expanding cone $\overline{\mathcal{C}'}$ for $H^{-1}$.

We may define local stable and unstable manifolds at any point $z$ where we have non-zero Lyapunov exponents. These are line segments aligned with the subspace $E_z^s$ as defined in Theorem \ref{thm:Oseledets-2d}, taking $f=H$ to find the stable direction, and $f=H^{-1}$ to find the unstable direction. The following lemma provides bounds on the gradients of these line segments.

\begin{lemma}
\label{lemma:alignment}
Given local unstable, stable manifolds $\gamma_u(z)$, $\gamma_s(z)$ at $z \in X'$, let $m_0$,$n_0$ be the smallest non-negative integers such that $H^{m_0}(z) \in H(A)$, $H^{-n_0}(z) \in H^{-1}(\mathfrak{a})$. Then
\begin{itemize}
    \item $H^{m_0}(\gamma_u(z))$ contains a segment $\gamma$ aligned with some vector $v \in \mathcal{C}$,
    \item $H^{-n_0}(\gamma_s(z))$ contains a segment $\gamma'$ aligned with some vector $v' \in \mathcal{C}'$.
\end{itemize}
\end{lemma}

\begin{proof}
We first note the link between the two minimal cones. Let $v_u(M_j)$, $v_s(M_j)$ be vector subspaces generated by the unstable and stable eigenvectors of some hyperbolic matrix $M_j$. Clearly $v_u(M_1) = v_s(M_1^{-1}) = v_s(\mathfrak{M}_1)$ and, in fact, we can always relate the stable, unstable eigenvectors of $M_j$ to the unstable, stable eigenvectors of $\mathfrak{M}_j$. For $j =2,3$ these are given by
\begin{equation}
\label{eq:eigenStable}
    v_s(M_j) = DH_1 v_u(\mathfrak{M}_j)
\end{equation}
and
\begin{equation}
\label{eq:eigenUnstable}
    v_u(M_j) = DH_1 v_s(\mathfrak{M}_j).
\end{equation}
To see this, note that in the $j=2$ case:
\begin{equation*}
    \begin{split}
        M_2^{-1} \cdot DH_1 \, v_u(\mathfrak{M}_2) & = DH_0^{-1}DH_1^{-1} \cdot DH_1 \,v_u(\mathfrak{M}_2) \\
        & = \left( DH_1 DH_1^{-1}\right) DH_0^{-1} \,v_u(\mathfrak{M}_2) \\
        & = DH_1 \mathfrak{M}_2\, v_u(\mathfrak{M}_2) \\
        & = c \, DH_1 \,v_u(\mathfrak{M}_2)
    \end{split}
\end{equation*}
for some $c$ with $|c|>1$. This implies $DH_1 v_u(\mathfrak{M}_2)$ is in the stable subspace of $M_2$, showing (\ref{eq:eigenStable}). The same argument applied to the right hand side of (\ref{eq:eigenUnstable}) yields $|c|<1$ as required. The case $j=3$ is analogous.

Now let $\gamma_u(z)$ be the local unstable manifold at some $z \in X'$. By the partition construction, $m_0$ is in $\{0,1,2\}$. Now $H^{m_0}(\gamma_u(z))$ is a piecewise linear curve, the union of at most 3 line segments $\gamma_j$. Since $z$ lies outside of the singularity set $S$, $H^{m_0}(z)$ lies in the interior of some $\gamma_j$, call it $\gamma$.

By definition, for any $\zeta,\zeta' \in \gamma_u(z)$
\begin{equation*}
    \mathrm{dist}(H^{-n}(\zeta),H^{-n}(\zeta')) \rightarrow 0
\end{equation*} 
as $n \rightarrow \infty$. By extension we have that
\begin{equation}
\label{eq:unstable}
    \mathrm{dist}(H^{-n}(\xi),H^{-n}(\xi')) \rightarrow 0
\end{equation} 
for any $\xi,\xi' \in \gamma$. 

This means that $H^{-1}(\gamma) \subset H^{-1}(\mathfrak{a})$ must be aligned with some vector in the cone region $\mathcal{C}_s$ bounded by $v_s(\mathfrak{M}_1)$ and $v_s(\mathfrak{M}_3)$, which includes $v_s(\mathfrak{M}_2)$\footnote{The argument for the $\varepsilon$-perturbation in section \ref{sec:eps} is analogous, but there the cone is bounded by $v_s(\mathfrak{M}_2)$ and $v_s(\mathfrak{M}_3)$}. For if it falls outside of this region, it will be pulled into the invariant expanding cone $\overline{\mathcal{C}'}$ for $H^{-1}$, which contradicts (\ref{eq:unstable}). Now if we apply $H$ to $H^{-1}(\gamma) \subset A$, $\gamma$ must align with a vector in $DH_1 \, \mathcal{C}_s$. By (\ref{eq:eigenUnstable}), this is precisely the minimal cone for $H$. The argument for local stable manifolds is analogous, instead using (\ref{eq:eigenStable}).

%This means that $H^{-1}(\gamma_u(z)) \subset H^{-1}(\mathfrak{a})$ must be aligned with some vector in the cone region $\mathcal{C}_s$ bounded by $v_s(\mathfrak{M}_1)$ and $v_s(\mathfrak{M}_3)$, which includes $v_s(\mathfrak{M}_2)$\footnote{The argument for the $\varepsilon$-perturbation is analogous, but there the cone is bounded by $v_s(\mathfrak{M}_2)$ and $v_s(\mathfrak{M}_3)$}. For if it falls outside of this region, it will be pulled into the invariant expanding (common) cone for $H^{-1}$, which contradicts (\ref{eq:unstable}). Now if we apply $H$ to $H^{-1}(\gamma_u(z)) \subset A$, $\gamma_u(z)$ must align with a vector in $DH_1 \, \mathcal{C}_s$. By (\ref{eq:eigenUnstable}), this is precisely the minimal cone for $H$. The argument for local stable manifolds is analogous, instead using (\ref{eq:eigenStable}).
\end{proof}

The main result of this section is the following.

\begin{prop}
\label{prop:etaM}
Condition \textbf{(M)} holds for $H$ when $0<\eta<\eta_1\approx 0.324$.
\end{prop}

\begin{figure}
    \centering

\begin{tikzpicture}

\node at (-3.6,0) {

    \begin{tikzpicture}[scale = 0.6]
    
    \draw (0,0) -- (0,7.5) -- (2.5,10) -- (10,10) -- (0,0);
    \draw (10,0) -- (10,7.5) -- (2.5,0) -- (10,0);
    
    \draw[fill=gray!50] (0,7.5) -- (2.5,7.5) -- (10,10) -- (7.5,10) -- (0,7.5);
    
    \node at (10.7,7.5) {$1-\eta$};
    \node at (7.5,10.3) {$1-\eta$};
    \draw[dashed] (0,0) -- (2.5,0);
    
    \draw[dashed] (0,7.5) -- (0,10) -- (2.5,10);
    
    \draw[dashed] (10,10) -- (10,7.5);
    
    \draw[fill=gray!20] (2.5,7.5) -- (7.5,7.5) -- (10,10) -- (2.5,7.5);
    \draw[fill=gray!20] (0,7.5) -- (7.5,10) -- (2.5,10) -- (0,7.5);
    
    \node at (2.5,10.3) {$\eta$};
    \node[scale=2] at (3,6) {$\mathfrak{a}_1$};
    \node[scale=2] at (7.5,2.5) {$\mathfrak{a}_1$};
    
    \node at (5,8.75) {$\mathfrak{a}_3$};
    \node at (3,9.25) {$\mathfrak{a}_2$};
    \node at (7,8.25) {$\mathfrak{a}_2$};

    \end{tikzpicture}
};
    
\node at (3.6,-0.15) {    
    \begin{tikzpicture}[scale = 0.6]
    
    \draw (0,0) rectangle (10,7.5);
    \draw [dashed] (0,7.5) -- (0,10) -- (10,10) -- (10,7.5);
    
    \draw[fill=gray!50] (0,0) -- (2.5,2.5) -- (10,7.5) -- (7.5,5) -- (0,0);
    \draw[fill=gray!20] (2.5,2.5) -- (10,7.5) -- (7.5,7.5) -- (2.5,2.5);
    \draw[fill=gray!20] (0,0) -- (7.5,5) -- (2.5,0) -- (0,0);
    
    \node at (2.5,-0.3) {$\eta$};
    \node at (-0.25,2.5) {$\eta$};
    \node at (10.7,7.5) {$1-\eta$};
    \node at (10.8,5) {$1-2\eta$};
    \node at (7.5,10.3) {$1-\eta$};
    
    \node[scale=2] at (3,6) {$A_1$};
    \node[scale=2] at (7.5,2.5) {$A_1$};
    
    \node at (5,3.75) {$A_3$};
    \node at (2.3,0.8) {$A_2$};
    \node at (7.5,6.5) {$A_2$};

    \end{tikzpicture}
};

\end{tikzpicture}

    \caption{Left: a partition of $\mathfrak{a}$ into three parts $\mathfrak{a}_i$, where $i$ is the return time of points in $\mathfrak{a}_i$ to $\mathfrak{a}$. Right: the equivalent plot for $A$, considering return times under $H^{-1}$. }
    \label{fig:subpartition}
\end{figure}
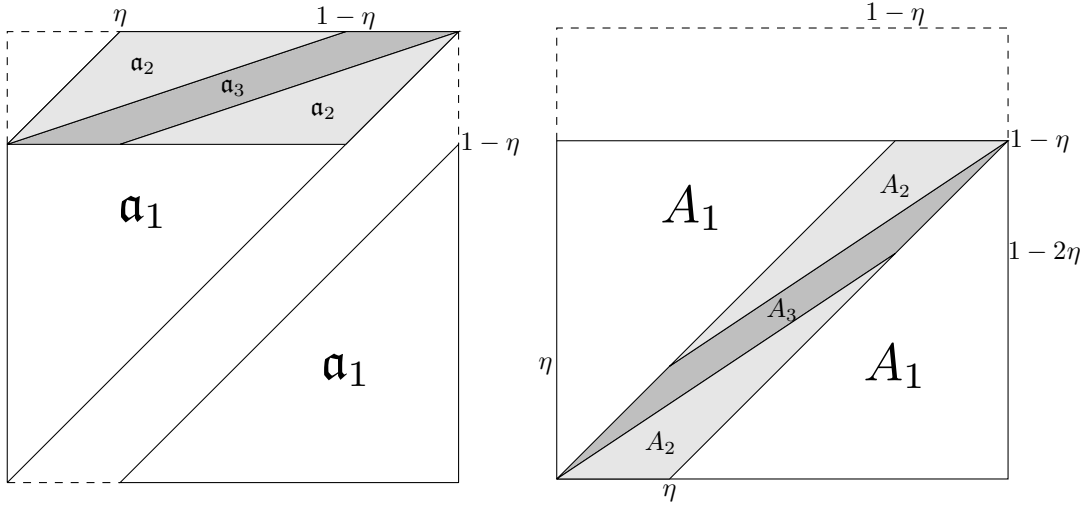

We will use the known behaviour of returns to $H(A)$ (resp. $H^{-1}(\mathfrak{a})$), and expansion during this return, to grow the images of local manifolds to the point where an intersection is certain in $A_1 = H(A) \cap H^{-1}(\mathfrak{a})$. This is a quadrilateral, shown in Figure \ref{fig:subpartition}. We call any line segment in $A_1$ which joins its upper and lower boundaries a $v$-segment. Similarly we call any line segment in $A_1$ which joins its left and right (sloping) boundaries a $h$-segment. Clearly $v$- and $h$-segments must always intersect. Given $z,z' \in X'$ our aim, then, is to find $m,n$ such that $H^m(\gamma_u(z))$ contains a $v$-segment and $H^{-n}(\gamma_s(z'))$ contains a $h$-segment.

The key issue we have to overcome in the growth stage is that while the images of the segments may grow exponentially in total length, the sign alternating property (as described in \citealp{cerbelli_continuous_2005}) means that they can repeatedly double back on themselves, meaning that their total diameter (be this in the $x$ or $y$ directions) does not necessarily grow. When considering unstable manifolds, we define the \emph{diameter} of a line segment $\Gamma$ as $\mathrm{diam}(\Gamma) = \nu \left( \{y \,|\, (x,y) \in \Gamma \} \right)$, where $\nu$ is the Lebesgue measure on $\mathbb{R}$. When considering stable manifolds, we instead define diameter using the projection to the $x$-axis.

We start with the method for growing unstable manifolds, partitioning $\mathfrak{a} = H(A)$ into three sets $\mathfrak{a}_i$, where the subscript $i$ is the return time of its elements to $\mathfrak{a}$. This is shown in Figure \ref{fig:subpartition}. We say that a line segment has \emph{non-simple intersection} with $\mathfrak{a}_i$ if its restriction to $\mathfrak{a}_i$ contains more than 1 connected component. The growth stage involves iteratively applying the following lemma.

\begin{lemma}
\label{lemma:etaUnstab}
Let $\Gamma_{p-1}$ be a line segment satisfying
\begin{enumerate}[label={(C\arabic*)}]
    \item $\Gamma_{p-1} \subset \mathfrak{a}$,
    \item $\Gamma_{p-1}$ is aligned with some vector in the minimal invariant cone $\mathcal{C}$ for $H$,
\end{enumerate}
and which has simple intersection with each of the $\mathfrak{a}_i$. There there exists a line segment $\Gamma_p$ satisfying (C1), (C2),
\begin{enumerate} [label={(C\arabic*)}]
    \setcounter{enumi}{2}
    \item $\Gamma_p \subset H^{i}(\Gamma_{p-1})$ for a chosen $i \in \{ 1,2,3\}$, and
    \item There exists $\delta>0$ such that $\mathrm{diam}(\Gamma_p) \geq (1+\delta)\,\mathrm{diam}(\Gamma_{p-1})$.
\end{enumerate}
\end{lemma}

% We will define a sequence of line segments $\Gamma_p$ such that for $0 \leq p \leq P$:
% \begin{enumerate}[label={(C\arabic*)}]
%     \item $\Gamma_p \subset \mathfrak{a}$,
%     \item $\Gamma_p$ is aligned with some vector in the minimal invariant cone for $H$,
% \end{enumerate}
% and for $1 \leq p \leq P$:
% \begin{enumerate} [label={(C\arabic*)}]
%     \setcounter{enumi}{2}
%     \item $\Gamma_p \subset H^{i}(\Gamma_{p-1})$ for a chosen $i \in \{ 1,2,3\}$,
%     \item There exists $\delta>0$ such that $\mathrm{diam}(\Gamma_p) \geq (1+\delta)\,\mathrm{diam}(\Gamma_{p-1})$,
% \end{enumerate}
% terminating in $\Gamma_P$ which has non-simple intersection with one of the $\mathfrak{a}_i$.

% Let $\gamma_u(z)$ be the local unstable manifold of some $z \in X'$. Let $n_0 \geq 0$ be the smallest integer such that $H^{n_0}(z) \in \mathfrak{a}$, then by Lemma \ref{lemma:alignment}, $H^{n_0}(\gamma_u(z))$ contains a segment $\Gamma_0$ in $\mathfrak{a}$, aligned with some vector in the invariant cone $\mathcal{C}$. 

\begin{proof}

The process of generating $\Gamma_p$ from $\Gamma_{p-1}$ is as follows. Based on the location of $\Gamma_{p-1}$ in $\mathfrak{a}$, we will restrict $\Gamma_{p-1}$ to one of the $\mathfrak{a}_i$ then map it forwards under $H^i$ to give $\Gamma_p$, satisfying (C3). By definition of the $\mathfrak{a}_i$, (C1) is satisfied. If $\Gamma_{p-1}$ is aligned with some $v \in \mathcal{C}$, $\Gamma_p$ is aligned with $M_i v$. By cone invariance, this is also in $\mathcal{C}$, so (C2) is satisfied. 

The expansion in diameter can be bounded from below by
\[ K_i(\eta) = \inf_{v \in \mathcal{C}} \frac{|| M_i v ||}{|| v ||} \]
where, again, we are using the $||\cdot||_\infty$ norm. Since we have already shown that the cone is expanding, if $\Gamma_{p-1}$ is entirely contained within some $\mathfrak{a}_i$ then taking $\Gamma_{p} = H^i(\Gamma_{p-1})$ ensures expansion in diameter. Where it becomes more interesting is when $\Gamma_{p}$ intersects multiple $\mathfrak{a}_i$. Looking at each of the $M_i$ across the invariant cone, at every parameter value $M_1$ has the smallest expansion on its eigenvector $v_u(M_1)$, $M_2$ and $M_3$ have the smallest expansion on the other cone boundary $v_u(M_3)$. Letting $\lambda_i$ be the magnitude of the unstable eigenvalue of $M_i$, $K_1$ and $K_3$ are given by 
\[K_1(\eta) = \lambda_1 (\eta) =  \frac{  3 - 2 \eta + \sqrt{5 - 4 \eta}}{2 \left(1-\eta \right)}
 \]
and
\[ K_3(\eta) = \lambda_3 (\eta) = \frac{ 3  - 9 \eta + 2\eta^{2} + \sqrt{- 36 \eta^{3} + 93 \eta^{2} - 54 \eta + 9} }{2\eta^{2}}.  \]
Next
\[ K_2(\eta) = \frac{2 \eta - 3}{1- \eta} \frac{1}{g_3^u(\eta)} + \frac{3- \eta}{\eta}, \]
calculated using the lower elements of $M_2$, the unit vector $\left(\frac{1}{g_3^u}, 1 \right)^T$, and the fact that $M_2$ reverses the orientation of vectors in the cone. 

Throughout, we assume that $\Gamma_{p-1}$ has simple intersection with each of the $\mathfrak{a}_i$. Suppose $\Gamma_{p-1}$ intersects $\mathfrak{a}_1$ and $\mathfrak{a}_2$, and write its restriction to these sets as $\Gamma^1$ and $\Gamma^2$ respectively. Since $K_1(\eta)$ and $K_2(\eta)$ are both greater than 2 for all $0<\eta< \frac{1}{3}$, and one of $\Gamma^1,\Gamma^2$ has diameter greater than or equal to $\frac{1}{2}$, we can restrict to that segment $\Gamma^i$ and expand under $H^i$ to establish that $\Gamma_{p}$ has larger diameter than $\Gamma_{p-1}$. Now suppose $\Gamma_{p-1}$ intersects $\mathfrak{a}_1$ and $\mathfrak{a}_3$. If the proportion of the diameter of $\Gamma_{p-1}$ in $\mathfrak{a}_1$ is greater than $\frac{1}{K_1(\eta)}$, we can simply expand from there. Otherwise $\Gamma^3$ has diameter greater than or equal to $1-\frac{1}{K_1(\eta)}$, and we can expand from $\mathfrak{a}_3$ provided that
\[ K_3(\eta) > \frac{1}{1-\frac{1}{K_1(\eta)}}. \]
The above is satisfied for approximately $\eta < 0.332$. The case where $\Gamma_{p-1}$ intersects $\mathfrak{a}_2$ and $\mathfrak{a}_3$ is similar and does not further restrict the parameter range.

Now suppose $\Gamma_{p-1}$ intersects $\mathfrak{a}_1$, $\mathfrak{a}_2$, and $\mathfrak{a}_3$. By the same argument as before, we require
\[ K_3(\eta)  > \frac{1}{1-\frac{1}{K_1(\eta)} - \frac{1}{K_2(\eta)} }. \]
Solving this numerically, the above inequality is satisfied for approximately $\eta < 0.327$. In any case, then, (C4) is satisfied.
\end{proof}

The method for growing the backwards images of local stable manifolds is entirely analogous. We divide up $A = H^{-1}(\mathfrak{a})$ into $A_1, A_2, A_3$ based on return time to $A$ under $H^{-1}$ (see Figure \ref{fig:subpartition}). The relevant hyperbolic matrices associated with the return map are $\mathfrak{M}_i$, which share an invariant, expanding cone $\mathcal{C}'$. We make minor adjustments to the (C) conditions to give:

\begin{lemma}
\label{lemma:etaStab}
Let $\Gamma_{p-1}$ be a line segment satisfying
\begin{enumerate}[label={(C\arabic*')}]
    \item $\Gamma_{p-1} \subset A$,
    \item $\Gamma_{p-1}$ is aligned with some vector in the minimal invariant cone $\mathcal{C}'$ for $H^{-1}$,
\end{enumerate}
and which has simple intersection with each of the $A_i$. There there exists a line segment $\Gamma_p$ satisfying (C1'), (C2'),
\begin{enumerate} [label={(C\arabic*')}]
    \setcounter{enumi}{2}
    \item $\Gamma_p \subset H^{-i}(\Gamma_{p-1})$ for a chosen $i \in \{ 1,2,3\}$,
    \item There exists $\delta>0$ such that $\mathrm{diam}(\Gamma_p) \geq (1+\delta)\,\mathrm{diam}(\Gamma_{p-1})$,
\end{enumerate}
where we measure the diameter of a line segment using its projection to the $x$-axis.
\end{lemma}

\begin{proof}

As before, define
\[ \mathcal{K}_i(\eta) = \inf_{v \in \mathcal{C}'} \frac{|| \mathfrak{M}_i v ||}{|| v ||}. \]
All of the $\mathfrak{M}_i$ see their minimum cone expansion on the cone boundary given by the unstable eigenvector of $\mathfrak{M}_3$. The key calculation we have to make is the parameter value $\eta_1$ such that
\begin{equation}
    \label{eq:K3K1K2}
    \mathcal{K}_3(\eta)  > \frac{1}{1-\frac{1}{\mathcal{K}_1(\eta)} - \frac{1}{\mathcal{K}_2(\eta)} }
\end{equation}
for $0<\eta < \eta_1$. We can solve numerically, giving $\eta_1 \approx 0.324$.
\end{proof}

Both of these lemmas hold, then, provided that $0<\eta < \eta_1$. They ensure the exponential growth in diameter of the segments $\Gamma_p$ up to some $\Gamma_P$ which has non-simple intersection with some $\mathfrak{a}_i$ (or $A_i$ for the stable case). At this point we will map directly into $v$- and $h$-segments.

% \subsubsection{Mapping into $v$-segments}

% Using the above, we can generate a sequence of line segments $\Gamma_{p}$ with exponentially increasing diameter. It follows that after a final number of steps, we reach a $\Gamma_P$ which no longer has simple intersection with all of the $\mathfrak{a}_i$ and we cannot keep growing by the above method. At this stage we will specify a specific number of iterations $k \in \{0,3,5\} $ such that $H^k(\Gamma_P)$ contains a $v$-segment.

\begin{lemma}
\label{lemma:etaVseg}
For any line segment $\Gamma_P \subset \mathfrak{a}$ which is aligned with a vector in $\mathcal{C}$ and has non-simple intersection with some $\mathfrak{a}_i$, $H^k(\Gamma_P)$ contains a $v$-segment for some $k \in \{0,3,5\} $.
\end{lemma}

\begin{proof}
All non-simple intersections give useful geometric information about $\Gamma_P$. Suppose it has non-simple intersection with $\mathfrak{a}_3$. Then as a connected straight line segment, it must \emph{traverse} $\mathfrak{a}_1$, that is, it connects the upper and lower boundaries of $\mathfrak{a}_1$, passing through $\mathfrak{a}_1$. By definition, this $\Gamma_P$ contains a $v$-segment. Now suppose $\Gamma_P$ has non-simple intersection with $\mathfrak{a}_2$. It follows that $\Gamma_P$ traverses $\mathfrak{a}_1$ ($v$-segment) or $\Gamma_P$ traverses $\mathfrak{a}_3$, connecting its sloping boundaries. This is case (I). We will show that any such segment contains a $v$-segment in its 5\textsuperscript{th} image. Finally assume that $\Gamma_P$ has non-simple intersection with $\mathfrak{a}_1$. It follows that we traverse $\mathfrak{a}_3$, case (I), or the restriction to $\mathfrak{a}_2$ is sufficiently large that its 3\textsuperscript{rd} image contains a $v$-segment, case (II).

We will start by showing case (I). Consider the quadrilateral $\mathcal{Q}_3 \subset \mathfrak{a}_3 $, defined by the four points $P_j$, shown in Figure \ref{fig:forwardsQ3}. Explicit coordinates for each of these points are given as part of the supplementary material. All of the points in the interior of $\mathcal{Q}_3$ share the same itinerary path under 5 iterations of $H$, $BCAAA$, so $H^5(\mathcal{Q}_3)$ is also a quadrilateral and any straight line segment contained within $\mathcal{Q}_3$ maps into a new straight line segment under $H^5$. It is clear that any $\Gamma_P$ which traverses $\mathfrak{a}_3$, joining its sloping boundaries, must also traverse $\mathcal{Q}_3$. The sloping boundaries of $\mathcal{Q}_3$ map into the upper and lower boundaries of $\mathfrak{a}_1$ under $H^5$, so if $\Gamma_P$ connects these sloping boundaries, $H^5(\Gamma_P)$ contains a $v$-segment.

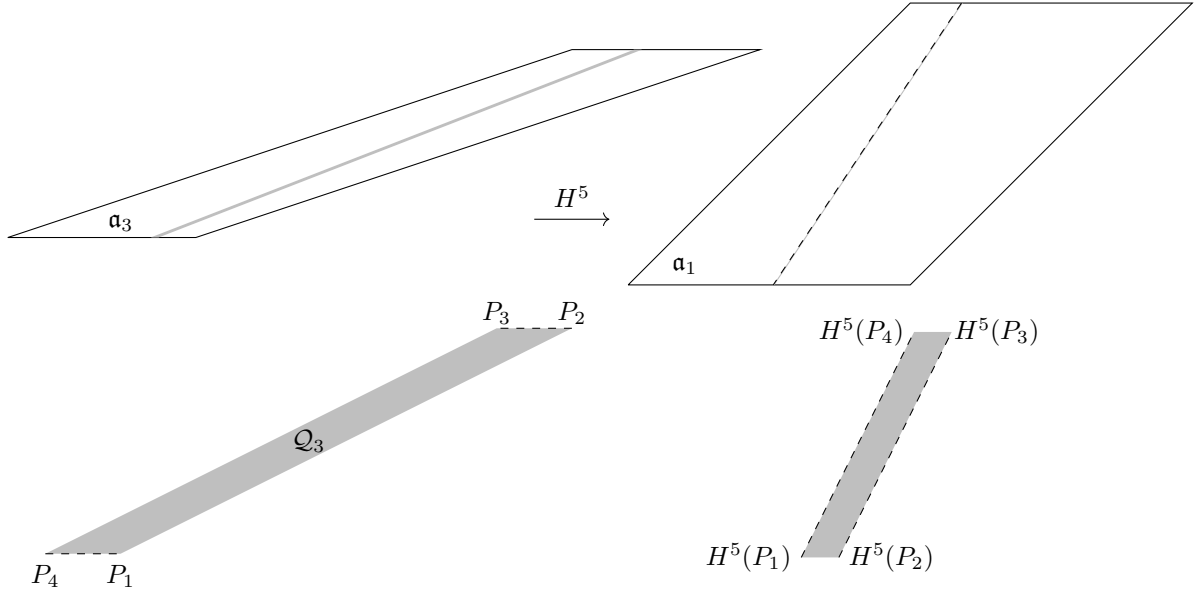
\begin{figure}
    \centering
\begin{tikzpicture}

\node at (-4,0) {

\begin{tikzpicture}

 \draw (0,7.5) -- (2.5,7.5) -- (10,10) -- (7.5,10) -- (0,7.5);
 
 \node at (1.5,7.7) {$\mathfrak{a}_3$};
 
 \fill [gray!50] (2.00339780449556,7.5) -- (1.89754312598014,7.5) -- (8.32592786199686,10) -- (8.43178254051228,10) --  (2.00339780449556,7.5);

\end{tikzpicture} };

\node at (-1.5,-0.7) {$H^5$};
\draw [->] (-2,-1) -- (-1,-1);

\node at (3,0) {
\begin{tikzpicture}[scale=0.5]

\draw (2.5,0) -- (10,0) -- (17.5,7.5) -- (10,7.5) -- (2.5,0);

\fill[gray!50] (6.35781495033978,0) -- (6.39309984317827,0) -- (1.39963408259278+10,7.5) -- (1.36434918975430+10,7.5) -- (6.35781495033978,0);

\draw[dashed] (1.36434918975430+10,7.5) -- (6.35781495033978,0);

\node at (4,0.5) {$\mathfrak{a}_1$};

\end{tikzpicture}

};

\node at (-5,-4) {
\begin{tikzpicture}
\fill [gray!50] (0,0) -- (1,0) -- (7,3) -- (6,3) -- (0,0);
\draw[dashed] (0,0) -- (1,0);
\draw[dashed] (7,3) -- (6,3);

\node at (3.5,1.5) {$\mathcal{Q}_3$};

\node at (1,-0.3) {$P_1$};

\node at (0,-0.3) {$P_4$};

\node at (7,3.2) {$P_2$};

\node at (6,3.2) {$P_3$};

\end{tikzpicture}

};

\node at (2.5,-4) {
\begin{tikzpicture}
\fill [gray!50] (0,0) -- (0.5,0) -- (2,3) -- (1.5,3) -- (0,0);
\draw[dashed] (0,0) -- (1.5,3);
\draw[dashed] (0.5,0) -- (2,3);

\node at (-0.7,0) {$H^5(P_1)$};

\node at (1.2,0) {$H^5(P_2)$};

\node at (2.6,3) {$H^5(P_3)$};

\node at (0.8,3) {$H^5(P_4)$};

\end{tikzpicture}

};

\end{tikzpicture}
    \caption{Case (I). A quadrilateral $\mathcal{Q}_3 \subset \mathfrak{a}_3$ and its image in $\mathfrak{a}_1$ under $H^5$. Any line segment $\Gamma$ which joins the sloping boundaries of $\mathfrak{a}_3$ will join the sloping boundaries of $\mathcal{Q}_3$, and hence $H^5(\Gamma \cap \mathcal{Q}_3)$ is a $v$-segment.}
    \label{fig:forwardsQ3}
\end{figure}

Case (II) can be argued similarly. We assume that $\Gamma_P$ has non-simple intersection with $\mathfrak{a}_1$ and that we do not traverse $\mathfrak{a}_3$ in such a way that we can argue as in case (I). We will concentrate first on the left portion of $\mathfrak{a}_2$; we shall soon see that the analysis for the right portion is analogous.

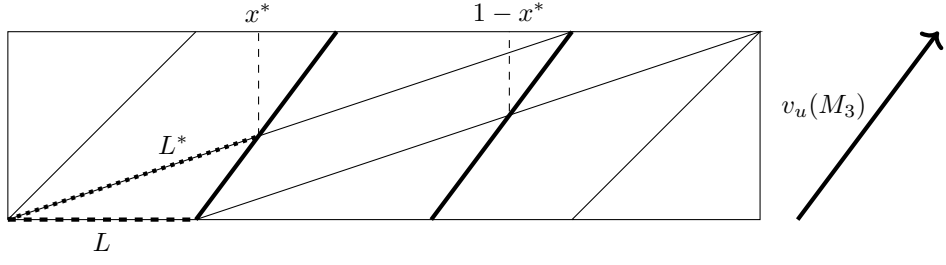
\begin{figure}
    \centering
    \begin{tikzpicture}
    
    \draw (0,7.5) -- (2.5,7.5) -- (10,10) -- (7.5,10);
    
    \draw[ultra thick] (2.5,7.5) -- (4.371238,10);
    \draw[ultra thick] (10-2.5,10) -- (10-4.371238,7.5);
    \draw (2.5,7.5) -- (7.5,7.5) -- (10,10) -- (10,7.5) -- (7.5,7.5);
    \draw  (0,7.5) -- (7.5,10) -- (2.5,10) -- (0,7.5);
    \draw (0,7.5) -- (0,10) -- (2.5,10);
    \draw[->, ultra thick] (8+2.5,7.5) -- (8+4.371238,10);
    \node at (10.85,9) {$v_u(M_3)$};
    
    \draw[dashed] (3.33311,8.61037) -- (3.33311,10);
    \draw[dashed] (10-3.33311,7.55 +10 - 8.61037) -- (10-3.33311,10);
    \node at (3.33311,10.25) {$x^*$};
    \node at (10-3.33311,10.25) {$1-x^*$};
    \draw[dashed, ultra thick] (0,7.5) -- (2.5,7.5);
    \node at (1.25,7.2) {$L$};
    \draw[dotted, ultra thick] (0,7.5) -- (3.33311,8.61037);
    \node at (2.2,8.5) {$L^*$};
    
    \end{tikzpicture}
    \caption{Geometry of line segments satisfying case (II). }
    \label{fig:xstar}
\end{figure}

Since we assume $\Gamma_P$ does not connect the sloping sides of $\mathfrak{a}_3$, it must intersect the $\mathfrak{a}_1$, $\mathfrak{a}_3$ boundary on $L$, shown in Figure \ref{fig:xstar}. The solid thick line shown is aligned with clockwise bound on the invariant cone, with gradient $g_3^u$. The intersection of $\Gamma_P$ with the $\mathfrak{a}_3$, $\mathfrak{a}_2$ boundary must lie in $L^*$, whose $x$-range is bounded above by $x^*$. 
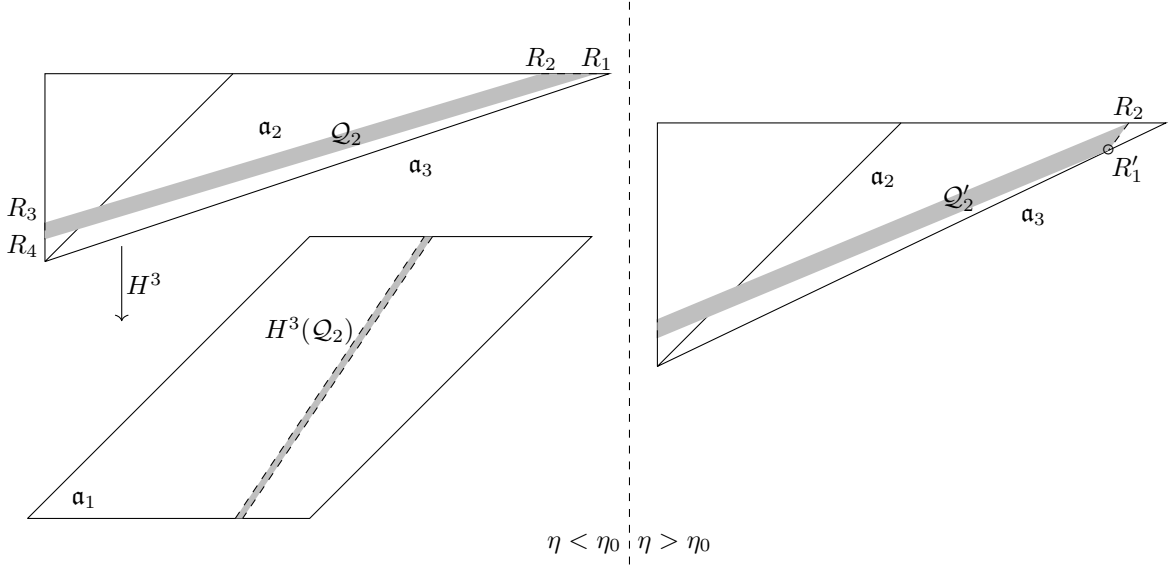
\begin{figure}[ht]
    \centering
\begin{tikzpicture}

\node at (-4,0) {

\begin{tikzpicture}

%\draw (0,7.5) -- (2.5,7.5) -- (10,10) -- (7.5,10) -- (0,7.5);

%\draw (2.5,7.5) -- (7.5,7.5) -- (10,10) -- (2.5,7.5);
\draw  (0,7.5) -- (7.5,10) -- (2.5,10) -- (0,7.5);
\draw (0,7.5) -- (0,10) -- (2.5,10);
%\draw (7.5,7.5) -- (10,7.5) -- (10,10);

    \node at (5,8.75) {$\mathfrak{a}_3$};
    \node at (3,9.25) {$\mathfrak{a}_2$};

\node at (7.33516483516484,10.2) {$R_1$};
\node at (6.59340659340659,10.2) {$R_2$};
\node at (-0.3,8.1980198019802) {$R_3$};
\node at (-0.3,7.70702970297030) {$R_4$};

\fill[gray!50] (0,8.01980198019802) -- (0,7.79702970297030) -- (7.33516483516484,10) -- (6.59340659340659,10) -- (0,8.01980198019802);

\node at (4,9.2) {$\mathcal{Q}_2$};

\draw[dashed] (0,8.01980198019802) -- (0,7.79702970297030);
\draw[dashed] (7.33516483516484,10) -- (6.59340659340659,10);

\end{tikzpicture} };

\node at (4,-1) {

\begin{tikzpicture}

%\draw (0,7.5) -- (2.5,7.5) -- (10,10) -- (7.5,10) -- (0,7.5);

%\draw (2.5,7.5) -- (7.5,7.5) -- (10,10) -- (2.5,7.5);
\draw  (0,10-3.23956002964573) -- (10-3.23956002964573,10) -- (3.23956002964573,10) -- (0,10-3.23956002964573);
\draw (0,10-3.23956002964573) -- (0,10) -- (3.23956002964573,10);
%\draw (7.5,7.5) -- (10,7.5) -- (10,10);

    \node at (5,8.75) {$\mathfrak{a}_3$};
    \node at (3,9.25) {$\mathfrak{a}_2$};

\fill[gray!50] (0,7.13563086647303) -- (0,7.38927641953645) -- (6.26611419495497,10) -- (5.99805335391734,9.63466916050191) -- (0,7.13563086647303);

\node at (4,9) {$\mathcal{Q}_2'$};

\node at (5.99805335391734,9.63466916050191) {$\circ$};

\node at (5.99805335391734+0.25,9.63466916050191-0.25) {$R_1'$};

\node at (6.26611419495497,10.2) {$R_2$};

\draw[dashed] (0,7.13563086647303) -- (0,7.38927641953645);
\draw[dashed] (6.26611419495497,10) -- (5.99805335391734,9.63466916050191);

\end{tikzpicture} };

\node at (-6.2,-1.75) {$H^3$};
\draw [->] (-6.5,-1.25) -- (-6.5,-2.25);

\draw[dashed] (0.25,2) -- (0.25,-5.5);
\node at (-0.35,-5.2) {$\eta <\eta_0$};
\node at (0.85,-5.2) {$\eta >\eta_0$};

\node at (-4,-3) {
\begin{tikzpicture}[scale=0.5]

\draw (2.5,0) -- (10,0) -- (17.5,7.5) -- (10,7.5) -- (2.5,0);

\fill[gray!50] (8.21782178217821,0) -- (8.02197802197803,0) -- (3.04945054945055+10,7.5) -- (3.26732673267326+10,7.5) -- (8.21782178217821,0);

\draw[dashed] (8.02197802197803,0) -- (3.04945054945055+10,7.5);
\draw[dashed] (3.26732673267326+10,7.5) -- (8.21782178217821,0);

\node at (10,5) {$H^{3}(\mathcal{Q}_2)$};

\node at (4,0.5) {$\mathfrak{a}_1$};

\end{tikzpicture}

};

\end{tikzpicture}
    \caption{Case (II) for $\eta$ either side of the critical value $\eta_0 = 1-\frac{1}{\sqrt{2}}$.}
    \label{fig:forwardsQ2}
\end{figure}

Let $\Gamma$ be the restriction to $\mathfrak{a}_2$. We will show that $\Gamma$, constrained by the $x^*$, intersects a quadrilateral whose image under $H^3$ stretches across $\mathfrak{a}_1$ in much the same way we saw in case (I). For $\eta \leq \eta_0 = 1-\frac{1}{\sqrt{2}} \approx 0.293$, such a quadrilateral $\mathcal{Q}_2$ exists and has all four corners on the lines $x=0$, $y=1$ (see left hand side of Figure \ref{fig:forwardsQ2}). Starting with the top-right and cycling anti-clockwise, these corners have coordinates
\[ R_1 = \left( \frac{- \eta^{3} + 7 \eta^{2} - 13 \eta + 7}{3 \eta^{2} - 10 \eta + 8} , 1 \right), \quad R_2 = \left( \frac{2 \left(2 \eta^{2} - 5 \eta + 3\right)}{3 \eta^{2} - 10 \eta + 8} , 1  \right), \]
\[ R_3 = \left( 0, \frac{5 \eta^{2} - 13 \eta + 8}{\eta^{2} - 7 \eta + 8}  \right), \text{ and } R_4 = \left( 0, \frac{- \eta^{3} + 7 \eta^{2} - 14 \eta + 8}{\eta^{2} - 7 \eta + 8}  \right). \]
Any line segment joining the $\mathfrak{a}_2,\mathfrak{a}_3$ boundary to the $\mathfrak{a}_2,\mathfrak{a}_1$ boundary must connect the parallel boundaries of $\mathcal{Q}_2$ and therefore maps into a $v$-segment. At the critical value $\eta = \eta_0$ the point $R_1$ lies on the rightmost corner of $\mathfrak{a}_2$, $(1-\eta,1)$. Now let $\eta > \eta_0$ and consider the quadrilateral $\mathcal{Q}_2'$ defined by the corners $R_2,R_3,R_4$, and
\begin{equation}
\label{eq:x'}
     R_1' = (x',y') = \left( \frac{- \eta^2 + 2 \eta -1}{\eta (2\eta-3)}, \frac{- 2 \eta^{2} + 6 \eta - 4}{2 \eta - 3} \right).
\end{equation}
This final corner also maps into $y=1-\eta$ under $H^3$, hence any line segment which joins the parallel sides of $\mathcal{Q}_2'$ maps into a $v$-segment. Certainly if $x^*(\eta)< x'(\eta)$ for $\eta_0<\eta<\eta_1$, then $\Gamma$ will connect the parallel sides of $\mathcal{Q}_2'$. First we solve line equations to give
\[ x^*(\eta) = \frac{\eta g_3^u(\eta)}{g_3^u(\eta) - \frac{\eta}{1-\eta} } \]
which is bounded from above by $x^*(\eta_1) \approx 0.5512$. Next by (\ref{eq:x'}),
\[ x'(\eta) = \frac{- \eta^2 + 2 \eta -1}{\eta (2\eta-3)} \]
which is bounded from below by $x'(\eta_1) \approx 0.5998$, establishing the result.

The case where $\Gamma$ traverses the other (right) part of $\mathfrak{a}_2$ is analogous. Note that we can transform one part of $\mathfrak{a}_2$ into the other by reflecting in the lines $y=1-\frac{\eta}{2}$ and $x=\frac{1}{2}$\footnote{Since the lines are orthogonal, $S_x \circ S_y = S_y \circ S_x$.}, written as $(S_x \circ S_y) (\mathfrak{a}_2) = \mathfrak{a}_2$. Now the images of $\mathcal{Q}_2$ and $\mathcal{Q}_2'$ under $S_x \circ S_y$ span across the right portion of $\mathfrak{a}_2$ in an analogous fashion to before and also map into $v$-segments under $H^3$. Making the same assumption as before, that case (II) holds but case (I) does not, we know that $\Gamma$ intersects the $\mathfrak{a}_2,\mathfrak{a}_3$ boundary at some point $(x,y)$ with $x>1-x^*$ (see Figure \ref{fig:xstar}). To ensure that $\Gamma$ connects the parallel sides of $(S_x \circ S_y) (\mathcal{Q}_2')$, it remains to check that the $x$-coordinate of $(S_x \circ S_y) (x',y')$, $1-x'$, is strictly less than $1-x^*$ across $\eta_0< \eta < \eta_1$. Indeed, $1-x'(\eta) < 1-x^*(\eta) $ follows from $x^*(\eta)< x'(\eta)$, established in the previous case.
\end{proof}

\begin{lemma}
\label{lemma:etaHseg}
For any line segment $\Gamma_P \subset A$ which is aligned with a vector in $\mathcal{C}'$ and has non-simple intersection with some $A_i$, $H^{-k}(\Gamma_P)$ contains a $h$-segment for some $k \in \{0,3,5\}$.
\end{lemma}

\begin{proof}

The argument is similar to the forwards-time case. A partition of $H^{-1}(\mathfrak{a}) = A$ by return time is shown in Figure \ref{fig:subpartition}. Case (I) assumes that $\Gamma$ connects the two $A_2,A_3$ boundaries through $A_3$, case (II) assumes that $\Gamma$ joins the two sloping boundaries of $A_1$ through $A_2 \cup A_3$, but that case (I) does not hold. We will show that in case (I) $H^{-5}(\Gamma)$ contains a $h$-segment, and in case (II) $H^{-3}(\Gamma)$ contains a $h$-segment. Starting with $\Gamma$ satisfying case (I), Figure \ref{fig:backwardsQuads} shows a quadrilateral $Q_3 \subset A_3$ with two short sides on the $A_1,A_3$ boundaries. It follows that $\Gamma$ must connect a segment which joins the longer sides of $Q_3$, through $Q_3$. The argument is now the same as in the forwards time analysis, all points in $Q_3$ share the same itinerary under 5 iterations of $H^{-1}$, $\mathfrak{b}\mathfrak{c}\mathfrak{a}\mathfrak{a}\mathfrak{a}$, so $H^{-5}(Q_3)$ is a quadrilateral in $A$. One can verify that it is wholly contained in $A_1 \subset A$ and that its longer sides map into its sloping boundaries (see right image in Figure \ref{fig:backwardsQuads}). $H^{-5}(\Gamma)$ then contains a segment which connects these two boundaries through $A_1$, that is, $H^{-5}(\Gamma)$ contains a $h$-segment. Explicit expressions for the corner coordinates of $Q_3$ and their images under $H^{-5}$ will be given as supplementary material.

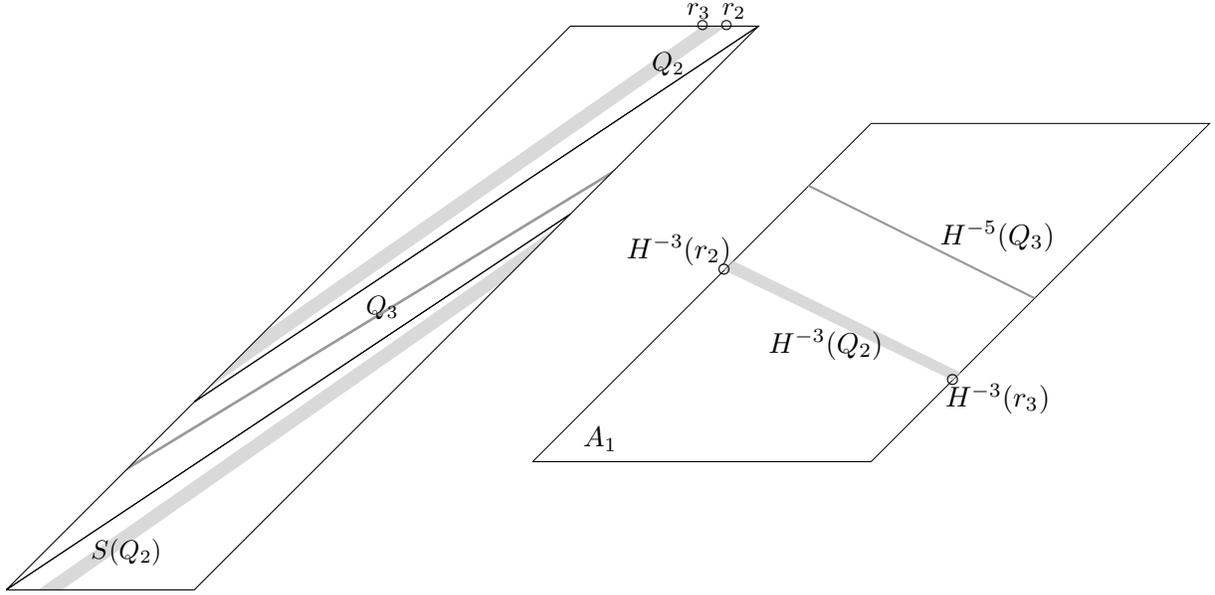
\begin{figure}
    \centering
    
    \begin{tikzpicture}
    
    \node at (-3.5,0) {

    \begin{tikzpicture}
    
    \fill[gray!80] (1.56821745948772,1.56821745948772) -- (7.99660219550444,7.99660219550444-2.5) -- (8.10245687401986,8.10245687401986-2.5) -- (1.67407213800314,1.67407213800314);
    
    \fill[gray!30] (2.66483516483516,2.66483516483516) -- (9.57547169811321,7.5) -- (9.25707547169811,7.5) -- (3.40659340659341,3.40659340659341);
    
    \fill[gray!30] (10-2.66483516483516,7.5-2.66483516483516) -- (10-9.57547169811321,0) -- (10-9.25707547169811,0) -- (10-3.40659340659341,7.5-3.40659340659341);
    
    \draw (0,0) -- (2.5,2.5) -- (10,7.5) -- (7.5,5) -- (0,0);
    \draw (2.5,2.5) -- (10,7.5) -- (7.5,7.5) -- (2.5,2.5);
    \draw (0,0) -- (7.5,5) -- (2.5,0) -- (0,0);
    
    \node at (9.67547169811321,7.7) {$r_2$};
    \node at (9.20707547169811,7.7) {$r_3$};
    
    \node at (9.57547169811321,7.5) {$\circ$};
    \node at (9.25707547169811,7.5) {$\circ$};
    
    \node at (5,3.75) {$Q_3$};
    
    \node at (8.8,7) {$Q_2$};
    
    \node at (1.6,0.5) {$S(Q_2)$};

    \end{tikzpicture}
    
    };
    
    \node[scale=0.6] at (3,0) {
    
    \begin{tikzpicture}
    \draw (2.5,0) -- (10,0) -- (17.5,7.5) -- (10,7.5) -- (2.5,0);
    
    \fill[gray!80] (10+3.60690015682176,3.60690015682176) -- (10+3.64218504966023,3.64218504966023) --  (8.63565081024572,8.63565081024572-2.5) -- (8.60036591740722,8.60036591740722-2.5);
    
    \fill[gray!30] (6.95054945054945,6.95054945054945-2.5) -- (6.74528301886793,6.74528301886793-2.5) -- (10+1.80424528301887 ,1.80424528301887 ) -- (10+1.97802197802198  , 1.97802197802198);
    
    \node[scale=1.8] at (6.74528301886793,6.74528301886793-2.5) {$\circ$};
    \node[scale=1.8] at (6.74528301886793-1,6.74528301886793-2.5 +0.45) {$H^{-3}(r_2)$};
    
    \node[scale=1.8] at (10+1.80424528301887 ,1.80424528301887 ) {$\circ$};
    \node[scale=1.8] at (10+1.80424528301887 +1 ,1.80424528301887 -0.4 ) {$H^{-3}(r_3)$};
    
    \node[scale=1.8] at (9,2.6) {$H^{-3}(Q_2)$};
    \node[scale=1.8] at (12.8,5) {$H^{-5}(Q_3)$};
    
    \node[scale=1.8] at (4,0.5) {$A_1$};
    
    \end{tikzpicture}   
    };
    
     \end{tikzpicture}
    
    \caption{Two quadrilaterals $Q_2 \subset A_2$ and $Q_3 \subset A_3$ which map into $A_1$ under $H^{-3}$ and $H^{-5}$ respectively. Their long boundaries map into the sloping boundaries of $A_1$, so segments $\Gamma$ which join these long boundaries map into $h$-segments. Case illustrated $\eta = \frac{1}{4}$.}
    \label{fig:backwardsQuads}
\end{figure}

Moving onto $\Gamma$ satisfying case (II) and first focusing on the upper portion of $A_2$, for $\eta \leq \eta_0$ we can follow the same argument, defining a quadrilateral $Q_2 \subset A_2$ with itinerary $\mathfrak{b}\mathfrak{a}\mathfrak{a}$ and $H^{-3}(Q_2) \subset A_1$ (see Figure \ref{fig:backwardsQuads}). Its long sides must be joined by $\Gamma$, and map into the boundary of $A_1$, so $H^{-3}(\Gamma)$ contains a $h$-segment. Starting with the bottom corner of $Q_2$ nearest the $A_2$, $A_3$ boundary and cycling anti-clockwise, label these points as $r_1, \dots , r_4$, which have coordinates
\[ r_1 = \left( \frac{\eta^{3} - 4 \eta^{2} + 3 \eta + 1}{3 \eta^{2} - 10 \eta + 8} , \frac{\eta^{3} - 4 \eta^{2} + 3 \eta + 1}{3 \eta^{2} - 10 \eta + 8}  \right), \quad r_2 = \left(  \frac{5 \eta^{3} - 20 \eta^{2} + 24 \eta - 8}{4 \eta^{3} - 18 \eta^{2} + 23 \eta - 8}  , 1-\eta \right), \]
\[ r_3 = \left( \frac{- \eta^{4} + 8 \eta^{3} - 23 \eta^{2} + 25 \eta - 8}{4 \eta^{3} - 18 \eta^{2} + 23 \eta - 8}, 1-\eta \right), \text{ and } r_4 = \left(  \frac{2 - \eta^{2}}{3 \eta^{2} - 10 \eta + 8} ,  \frac{2 - \eta^{2}}{3 \eta^{2} - 10 \eta + 8}  \right). \]
For $\eta > \eta_0$ we consider the quadrilateral $Q_2'$ with corners $r_2,r_3,r_4$ and
\[ r_1' = \left(   \frac{3 \eta^{2} - 5 \eta + 1}{\eta \left(2 \eta - 3\right)} ,  \frac{- 2 \eta^{3} + 7 \eta^{2} - 6 \eta + 1}{\eta \left(2 \eta - 3\right)} \right).\]
This is shown in Figure \ref{fig:q2dash}, with the $x$-coordinate of $r_1'$ highlighted as $x'(\eta)$. Like in the forwards-time case, we need to check that $x'(\eta)$ is not so far along the $A_2,A_3$ boundary that any $\Gamma$ satisfying case (II) does not connect the parallel sides of $Q_2'$. Letting $\mathfrak{g}_3^u(\eta)$ be the gradient of the anti-clockwise invariant cone boundary for $H^{-1}$, this amounts to showing that $x'(\eta) < x^*(\eta)$ where $(x^*,y^*)$ lies on the intersection of the lines
\[ y = \eta + \frac{1-2\eta}{1-\eta} (x-\eta) \]
(the $A_2,A_3$ boundary) and
\[ y = 1-2\eta + \mathfrak{g}_3^u(\eta) (x-1+\eta), \]
shown as the solid bold line in Figure \ref{fig:q2dash}. Solving for $x$ gives
\[ x^*(\eta) = \frac{\eta^2+3\eta-1+\mathfrak{g}_3^u(\eta)(1-\eta)^2}{\mathfrak{g}_3^u(\eta)(1-\eta) -1+2\eta}. \]
One can now verify that $x'(\eta) < x'(\eta_1) < x^*(\eta_1) < x^*(\eta)$ for all $\eta_0 < \eta < \eta_1$, establishing the result. To conclude case (II) we must extend the analysis to the other portion of $A_2$. This process is entirely analogous to the forwards time case, taking reflections in $x=\frac{1}{2}$ and $y=\frac{1}{2}-\frac{\eta}{2}$. An example is shown in Figure \ref{fig:backwardsQuads}, with the image of $Q_2$ under these reflections shown as $S(Q_2)$.

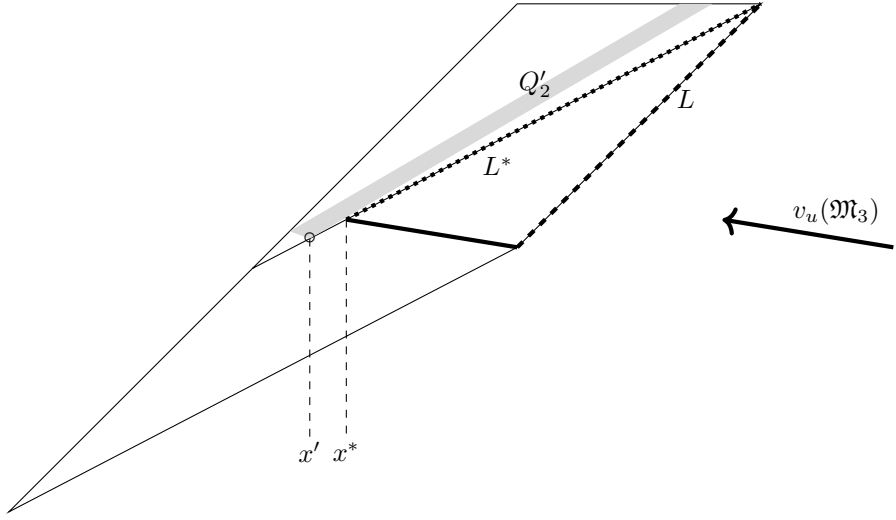
\begin{figure}
    \centering
    \begin{tikzpicture}
    \draw (3.23956002964573,3.23956002964573) -- (10,10-3.23956002964573) --(10-3.23956002964573,10-2*3.23956002964573) -- (0,0) -- (10-3.23956002964573,10-3.23956002964573) -- (10,10-3.23956002964573);
    
    \fill[gray!30] (4.00194664608266,3.63661580658458) -- (9.35684329038124,6.760439970354271) -- (8.92204105769041,6.760439970354271) -- (3.73388580504503,3.73388580504503);
    
    \node at (4.00194664608266,3.63661580658458) {$\circ$};
    
    \draw[ultra thick] (10-3.23956002964573, 10-2*3.23956002964573) -- (4.487799796343105,3.8896512034497914);
    
    \draw[dashed, ultra thick] (10-3.23956002964573, 10-2*3.23956002964573) -- (10, 10-3.23956002964573);
    \node at (9,5.5) {$L$};
    
    \draw[dotted, ultra thick] (4.487799796343105,3.8896512034497914) -- (10, 10-3.23956002964573);
    \node at (6.5,4.6) {$L^*$};
    \node at (7,5.7) {$Q_2'$};
    \draw [dashed] (4.00194664608266,3.63661580658458) -- (4.00194664608266,1);
    \draw [dashed] (4.487799796343105,3.8896512034497914) -- (4.487799796343105,1);
    
    \node at (4.00194664608266,0.8) {$x'$};
    \node at (4.487799796343105,0.8) {$x^*$};
    
    \draw [->, ultra thick] (5+10-3.23956002964573, 10-2*3.23956002964573) -- (5+4.487799796343105,3.8896512034497914);
    
    \node at (11,4) {$v_u(\mathfrak{M}_3)$};
    
    \end{tikzpicture}
    
    \caption{Case (II) for $\eta > \eta_0$. Any $\Gamma$ satisfying case (II) must intersect the $A_1 A_3$ boundary on $L$ and the $A_2 A_3$ boundary on $L^*$. This gives a lower bound on $x^*$ on the $x$-coordinate of this intersection so that if $x^* > x'$, then $\Gamma$ joins the parallel sides of $Q_2'$.} 
    \label{fig:q2dash}
\end{figure}
\end{proof}

We are now ready to prove Proposition \ref{prop:etaM}.

% To conclude, at almost every $z \in \tor$ we have found $m_0 \in \{ 0,1,2 \}$ such that $H^{m_0}(z) \in H(A)$, $m_1$ such that $H^{m_0 +m_1} (\gamma_u(z))$ contains a segment having non-simple intersection with one of the $\mathfrak{a}_i$, then $m_2 \in \{ 0,3,5 \}$ such that $H^{m_0 +m_1 +m_2}(\gamma_u(z))$ contains a $v$-segment. Similarly we can find $n = n_0 + n_1 + n_2$ such that $H^{-n}(\gamma_s(z'))$ contains a $h$-segment, for almost every $z' \in \tor$. Let $m = m_0+m_1+m_2$, it follows that $H^{m}(\gamma_u(z)) \cap H^{-n}(\gamma_s(z')) \neq \varnothing$, establishing \textbf{(M)}.

\begin{proof}[Proof of Proposition \ref{prop:etaM}]

Let $\gamma_u(z)$ be the local unstable manifold at some $z \in X'$. Let $m_0 \geq 0$ be the smallest integer such that $H^{m_0}(z) \in \mathfrak{a}$. Then by Lemma \ref{lemma:alignment}, $H^{m_0}(\gamma_u(z))$ contains a segment $\Gamma_0$ in $\mathfrak{a}$, aligned with some vector in the invariant cone $\mathcal{C}$. We can then iteratively apply Lemma \ref{lemma:etaUnstab} to generate a sequence of line segments with exponentially increasing diameter $(\Gamma_p)_{0\leq p \leq P}$ with each $\Gamma_p \subset H^{m_0 + m_p}(\gamma_u(z))$ for some $m_p>0$. Since the sequence has exponentially increasing diameter, after some finite number of steps $P$, the line segment $\Gamma_P$ must have non-simple intersection with one of the $a_i$. Lemma \ref{lemma:etaVseg} then tells us that $H^{k}(\Gamma_P)$ contains a $v$-segment for some $k \in \{0,3,5\}$. It follows that $H^{m}(\gamma_u(z))$ contains a $v$-segment where $m=m_0+m_P+k$. Similarly given $z' \in X'$, we can apply Lemmas \ref{lemma:alignment}, \ref{lemma:etaStab}, and \ref{lemma:etaHseg} to find $n$ such that $H^{-n}(\gamma_s(z'))$ contains a $h$-segment. Since $z$ and $z'$ were arbitrary, condition \textbf{(M)} holds.

\end{proof}

This establishes $H$ as ergodic over $0 <\eta < \eta_1$. Stronger mixing properties can now be easily shown.

%% file: sections/etaPert/repeatedManifold.tex
\begin{prop}
\label{prop:etaMR}
Condition \textbf{(MR)} holds for $H$ when $0<\eta<\eta_1\approx 0.324$.
\end{prop}

\begin{figure}
    \centering
    
    \begin{tikzpicture}
    
    \node[scale=0.5] at (-4.5,0) {
    
    \begin{tikzpicture}

    \fill[gray!20] (7,4.5) -- (12.25,2.25) -- (10,0) -- (4.75,2.25);
    %\fill[gray!20] (10,7.5) -- (15.25,5.25) -- (10,0) -- (4.75,2.25);
    \draw (2.5,0) -- (10,0) -- (17.5,7.5) -- (10,7.5) -- (2.5,0);
    
    \draw[ultra thick,dashed] (7,4.5) -- (12.25,2.25);
    \draw[ultra thick,dashed] (10,0) -- (4.75,2.25);
    
    \node[scale=1.8] at (8.5,2.5) {$Q^+$};
    
    \node[scale=1.8] at (4.25,2.25) {$c_1^+$};

    \end{tikzpicture}
    
    };
    
   \node[scale=0.5] at (0,-0.22) {
    
    \begin{tikzpicture}
    
    \fill[gray!40] (10,7.5) -- (12.25,7.5) -- (7,0) -- (4.75,0);
    
    \draw (2.5,0) -- (10,0) -- (17.5,7.5) -- (10,7.5) -- (2.5,0);
    
    \draw[ultra thick,dashed] (12.25,7.5) -- (7,0);
    \draw[ultra thick,dashed] (10,7.5)-- (4.75,0);
    
    \node[scale=1.8] at (8.5,3.75) {$Q^-$};
    
    \node[scale=1.8] at (4.75,-0.35) {$c_1^-$};
    
    \end{tikzpicture}
    
    };

    \node[scale=0.5] at (4.5,0) {
    
    \begin{tikzpicture}

    \fill[gray!40,opacity=0.5] (13,7.5) -- (15.25,7.5) -- (10,0) -- (7.75,0);
    
    \fill[gray!80,opacity=0.5] (10,7.5) -- (15.25,5.25) -- (13,3) -- (7.75,5.25);
    
    \draw (2.5,0) -- (10,0) -- (17.5,7.5) -- (10,7.5) -- (2.5,0);
    
        \draw[ultra thick,dashed] (13,7.5) -- (15.25,7.5);
    \draw[ultra thick,dashed] (10,0) -- (7.75,0);
    
        \draw[ultra thick,dashed] (10,7.5) -- (7.75,5.25);
    \draw[ultra thick,dashed] (15.25,5.25) -- (13,3) ;
    
    \node[scale=1.8] at (10,6) {$H^{-1}(Q^-)$};
    
    \node[scale=1.8] at (10,2.75) {$H(Q^+)$};
        
    \end{tikzpicture}
    
    };

    \end{tikzpicture}

    \caption{Two quadrilaterals $Q^+,Q^-$ in $A_1$ which map into $A_1$ under $H$ and $H^{-1}$ respectively. Any $v$-segment must join the dotted sides of $Q^+$, hence maps into another $v$-segment. Similar for $h$-segments and $Q^-$.}
    \label{fig:repeatedInt}
\end{figure}
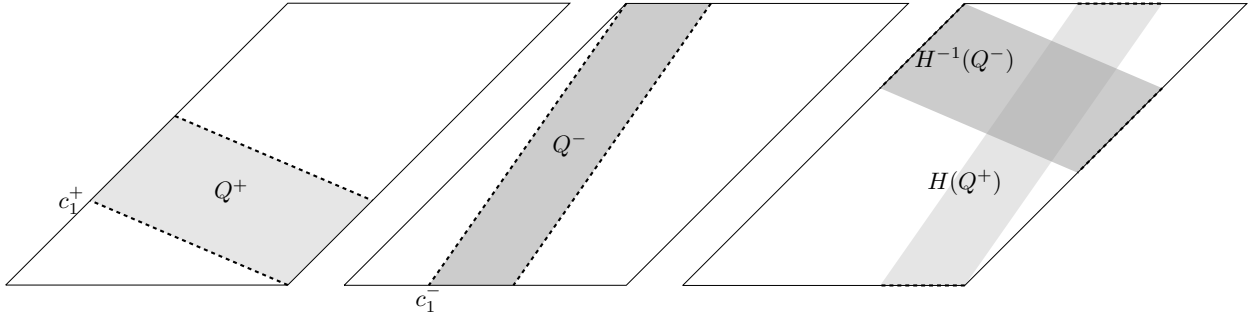

\begin{proof}
To establish \textbf{(MR)} it is sufficient to show that the image of a $v$-segment under $H$ contains a $v$-segment, and the image of a $h$-segment under $H^{-1}$ contains a $h$-segment. We can approach this is same way as before, defining quadrilaterals which these segments must traverse and looking at their images. Define the quadrilateral $Q^+$ by the corners (starting from the leftmost and cycling anti-clockwise)
\[ c_1^+ = \left( \frac{1+\eta-\eta^2}{3-2\eta}, \frac{(1-\eta)^2}{3-2\eta} \right) , c_2^+ = \left( 0,0 \right),  \]
\[ c_3^+ = \left( \frac{(1-\eta)^2}{3-2\eta} , \frac{(1-\eta)^2}{3-2\eta} \right) ,\text{ and }c_4^+ = \left( \frac{2-\eta}{3-2\eta} \frac{2(1-\eta)^2}{3-2\eta} \right).  \]
This is shown in the first diagram in Figure \ref{fig:repeatedInt}, note that the we have shifted the domain horizontally to more easily see $A_1$ as a quadrilateral. Any $v$-segment must join the dotted sides of $Q^+$, which map into the upper and lower boundaries of $A_1$, so $v$-segments map into $v$-segments. We can similarly define the quadrilateral $Q^-$ by the corners (starting from the leftmost and cycling anti-clockwise)
\[ c_1^- = \left( \frac{1+\eta-\eta^2}{3-2\eta}, 0\right), \, c_2^- = \left( \frac{2-\eta}{3-2\eta}, 0\right),\, c_3^- = \left( \frac{(1-\eta)^2}{3-2\eta} , 1-\eta \right) ,\text{ and }c_4^- = (0,1-\eta).  \]
Again, $h$-segments must connect the dotted sides of $Q^-$, which map into the sloping boundaries of $A_1$, hence $h$-segments map into $h$-segments.
\end{proof}

%% file: sections/epsPert/epsLyp.tex
Let $\varepsilon = \frac{1}{2}-\eta$. Our method for establishing non-zero Lyapunov exponents almost everywhere for $H$ as an $\varepsilon$-perturbation is essentially the same as in section \ref{sec:etaLyp}.
\begin{prop}
\label{prop:epsLyp}
We have non-zero Lyapunov exponents $\chi(z,v) \neq 0$ for almost every $z\in \tor $, $v\neq 0$, when $0<\varepsilon<\varepsilon_1 \approx 0.0931 $.
\end{prop}

\begin{proof}

The partition and possible itinerary paths $I_j$ around the partition are the same as before. Define the corresponding matrices $M_j$ using the derivative matrices
\[
DH_0 = \begin{pmatrix} 1 & \frac{-2}{1-2\varepsilon} \\ 1 & \frac{-1-2\varepsilon}{1-2\varepsilon} \end{pmatrix} \text{ and }
DH_1 = \begin{pmatrix} 1 & \frac{2}{1+2\varepsilon} \\ 1 & \frac{3+2\varepsilon}{1+2\varepsilon} \end{pmatrix}.
\]
Again, $M_3$ is the matrix which dictates our parameter range. It is hyperbolic for $\varepsilon< \varepsilon_1$, where $\varepsilon_1 = \frac{\sqrt{33}-5}{8} \approx 0.0931$. $M_2$ is hyperbolic for $\varepsilon$ strictly greater than 0.

Following the same argument as in section \ref{sec:etaLyp}, it remains to define an invariant cone and show that it is expanding. Defining $g_j^u$ and $g_j^s$ as before, one can verify that
\[ g_3^s(\varepsilon) < g_1^s(\varepsilon) < g_2^s(\varepsilon) < g_2^u(\varepsilon) < g_1^u(\varepsilon)< g_3^u(\varepsilon) \]
for $0< \varepsilon < \frac{1}{\sqrt{3}}- \frac{1}{2} \approx 0.0774$, and
\[ g_1^s(\varepsilon) < g_2^s(\varepsilon) < g_2^u(\varepsilon) < g_1^u(\varepsilon)< g_3^u(\varepsilon) < g_3^s(\varepsilon)  \]
for $\frac{1}{\sqrt{3}}- \frac{1}{2}< \varepsilon < \varepsilon_1$. Hence the cone $\mathcal{C}$, bounded by and including the unstable eigenvectors of $M_2$ and $M_3$, is the minimal invariant cone. The common cone $\overline{\mathcal{C}}$ is then defined as the open region bounded by the unstable eigenvector of $M_2$ at $\varepsilon=0$ and the unstable eigenvector of $M_3$ at $\varepsilon=\varepsilon_1$. Under the $||\cdot||_\infty$ norm, these are the unit vectors $(1,1)^T$ and $\left(\frac{\sqrt{33}-3}{6},1 \right)^T$ respectively. One can show that
\begin{multicols}{2}
\begin{itemize}
    \item $ ||M_1(1,1)^T||  > \frac{\sqrt{33}+9}{4}$
    \item $ ||M_2(1,1)^T|| > 1 $
    \item $ ||M_3(1,1)^T||  > \frac{9+\sqrt{33}}{6} $
    \item $ ||M_1\left(\frac{\sqrt{33}-3}{6},1 \right)^T|| >  \frac{9+5 \sqrt{33}}{12} $
    \item $ ||M_2\left(\frac{\sqrt{33}-3}{6},1 \right)^T||  > 7 - \frac{2 \sqrt{33}}{3} $
    \item $ ||M_3\left(\frac{\sqrt{33}-3}{6},1 \right)^T|| > 1 $
\end{itemize}
\end{multicols}
\noindent for all $\varepsilon$ in our range, so that our cone is expanding.
\end{proof}

This establishes non-uniform hyperbolicity. As before, the next section shows ergodicity.

%% file: sections/epsPert/epsMan.tex
\begin{prop}
\label{prop:epsM}
Condition \textbf{(M)} holds for $H$ over $\varepsilon_0<\varepsilon \leq \varepsilon_2$, where $\varepsilon_0 \approx 0.00925$ and $\varepsilon_2 \approx 0.0850$.
\end{prop}

The overall method for establishing \textbf{(M)} is unchanged. The key constructions are the partitions of $H(A)$ and $H^{-1}(\mathfrak{a})$ given in section \ref{sec:etaMan}, and the invariant cones $\mathcal{C}$ for $H$ (given above) and $\mathcal{C}'$ for $H^{-1}$. Defining the $\mathfrak{M}_j$ as before, $\mathcal{C}'$ is defined at each $\varepsilon$ as the cone bounded by (and including) the unstable eigenvectors of $\mathfrak{M}_2$ and $\mathfrak{M}_3$, i.e. the non-zero vectors with gradient $\mathfrak{g}_3^u<\mathfrak{g}<\mathfrak{g}_2^u$. One can show (by the same method as before) that $\mathcal{C}'$ is invariant and expanding.

For the sake of brevity, we will only describe the process of growing the backwards images of local stable manifolds. The process for unstable manifolds is entirely analogous and, due to $\mathcal{C}$ covering a smaller gradient range than $\mathcal{C}'$, results in less stringent bounds on the parameter range.

\begin{figure}
    \centering
    \subfigure[][]{%
\label{fig:Apartition}%
%\hspace{-3em}  
    \begin{tikzpicture}[scale = 0.85]
    \draw (0,0) rectangle (10,5.5);
    \draw (0,5.5) -- (0,6);
    \draw (10,5.5) -- (10,6);
    \draw[fill=gray!20] (0,0) -- (5.5,5.5) -- (10,5.5) -- (4.5,0) -- (0,0);
    \draw (0,5.5) -- (5.5,5.5);
   \node at (10.6,5.5) {$\frac{1}{2}+\varepsilon$};
   \draw[fill=gray!40] (0,0) -- (5.5,1) -- (4.5,0) -- (0,0);

   \draw[fill=gray!40] (10,5.5) -- (4.5,4.5) -- (5.5,5.5) --(10,5.5);
 %  \draw[dashed] (4.5,4.5) -- (0,4.5);
  % \draw[dashed] (5.5,1) -- (10,1);
   \node at (4.5,-0.4) {$\frac{1}{2}-\varepsilon$};
   \node at (10.4,1) {$2\varepsilon$};
      \node at (-0.55,4.5) {$\frac{1}{2}-\varepsilon$};
    \node at (5.5,5.8) {$\frac{1}{2}+\varepsilon$};
    \node at (5,2.5) {$A_3$};
    \node at (1.5,3) {$A_1$};
    \node at (8.5,2) {$A_1$};
    \node at (4,0.4) {$A_2$};
    \node at (6,5.1) {$A_2$};
%\draw[dashed] (5.5,1) -- (5.5,5.5);
 %  \draw[dashed] (4.5,4.5) -- (4.5,0);
   \end{tikzpicture}
   }%
   \vspace{1em}
\subfigure[][]{%
\label{fig:greenRed}%
 \begin{tikzpicture}[scale = 1.8]
   % \draw[fill=gray!20] (0,0) -- (5.5,5.5) -- (10,5.5) -- (4.5,0) -- (0,0);

   \draw[fill=gray!40] (0,0) -- (5.5,1) -- (4.5,0) -- (0,0);

    \draw[fill=gray] (0,0) -- (363/122,0) -- (110/21,110/21-4.5) -- (5.5,1)  -- (0,0);

    \draw[fill=gray!20] (6.2,1.7) -- (5.5,1) -- (0,0) -- (1.7,1.7);
    
    %\draw[dashed] (110/21,110/21-4.5) -- (6.2,110/21-4.5);
    
    \node at (110/21+0.2,110/21-4.55) {$Q_2$};
    
    \draw (0,0) -- (0,1.7);
    \draw (4.5,0) -- (6.2,0);
    
    \node at (363/122,-0.2) {$Q_1$};
     \node at (3.1,1.2) {$A_3$};
    % \node at (8.5,2) {$A_1$};
     \node[white] at (3.5,0.4) {$A_4$};
     \node at (4.25,0.2) {$A_5$};
    % \node at (6,5.1) {$A_2$};
   \end{tikzpicture}
    }
    \caption{Part (a) gives partition of $A$ based on return time to $A$ under iterations of $H^{-1}$. Part (b) shows a subdivision $A_4 \cup A_5 = A_2$, with the boundary between these sets defined as the segment joining the points $Q_1,Q_2$. Case illustrated $\varepsilon=0.05$.}
    \label{fig:A2partition}
\end{figure}
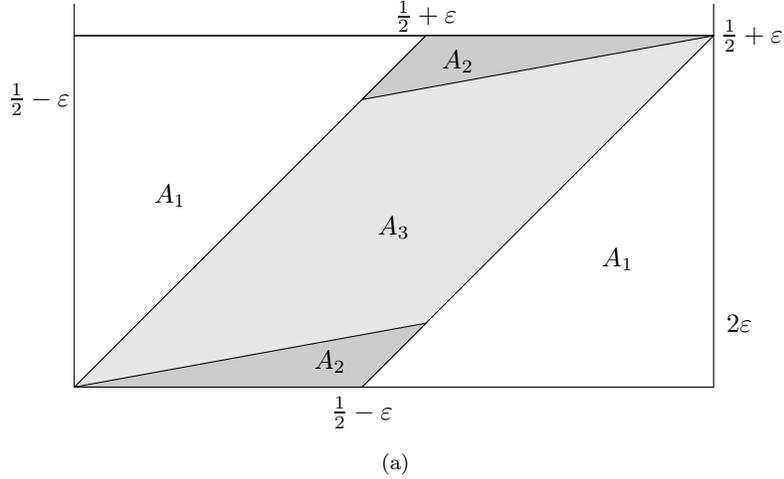
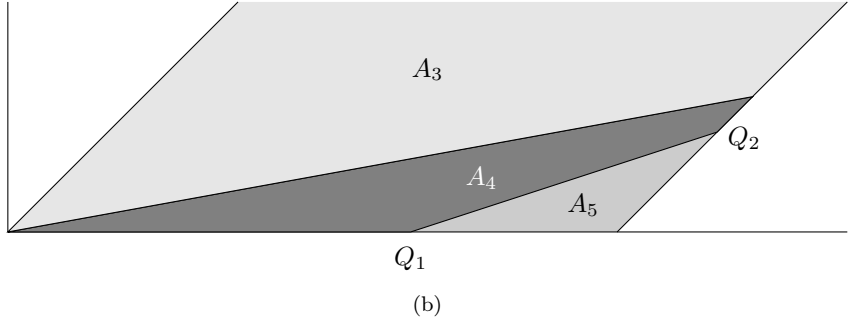

While for the $\eta$-perturbation the growth stage was relatively straightforward and the $h$-segment mappings more involved, the opposite is true for the $\varepsilon$-perturbation. If we were to follow the same method as before, reducing the parameter range to satisfy equations like (\ref{eq:K3K1K2}), we would be left with just a fragment of the parameter range. Our way around this necessitates growing piecewise linear curves rather than line segments. To ensure that we can find the diameter of a curve by summing the diameters of its constituent line segments, we require that a curve \emph{does not double back on itself}, that is, the projection to the $x$-axis is injective. The lemma for the growth stage is as follows:

\begin{lemma}
\label{lemma:epsStab}
Let $\Gamma_{p-1}$ be a piecewise linear curve satisfying
\begin{enumerate}[label={(C\arabic*')}]
    \setcounter{enumi}{-1}
    \item $\Gamma_{p-1}$ does not double back on itself,
    \item $\Gamma_{p-1} \subset A$,
    \item Each line segment in $\Gamma_{p-1}$ is aligned with some vector in the minimal invariant cone $\mathcal{C}'$ for $H^{-1}$,
\end{enumerate}
and which has simple intersection with each of the $A_i$. There there exists a piecewise linear curve $\Gamma_p$ satisfying (C0'), (C1'), (C2'),
\begin{enumerate} [label={(C\arabic*')}]
    \setcounter{enumi}{2}
    \item $\Gamma_p \subset H^{-i}(\Gamma_{p-1})$ for a chosen $i \in \{ 1,2,3\}$,
    \item There exists $\delta>0$ such that $\mathrm{diam}(\Gamma_p) \geq (1+\delta)\,\mathrm{diam}(\Gamma_{p-1})$,
\end{enumerate}
where we measure the diameter of a curve using its projection to the $x$-axis.
\end{lemma}

% \begin{enumerate}[label={(C\arabic*)}]
%     \item $\Gamma_p \subset A$,
%     \item Each line segment within $\Gamma_p$ is aligned with some vector in $\mathcal{C}'$,
%     \item $\Gamma_p \subset H^{-j}(\Gamma_{p-1})$ for a chosen $j \in \{ 1,2,3\}$,
%     \item There exists $\delta>0$ such that $\mathrm{diam}(\Gamma_p) \geq (1+\delta)\,\mathrm{diam}(\Gamma_{p-1})$, measuring diameter using the projection to the $x$-axis.
%     \item $\Gamma_p$ does not \emph{double back on itself}, i.e. its projection to the $x$-axis is injective.
%     \end{enumerate}

\begin{proof}

Figure \ref{fig:Apartition} shows the return time partition of $A = H^{-1}(\mathfrak{a})$ under $H^{-1}$. Define $\mathcal{K}_j(\varepsilon)$ for $j=1,2,3$ as before. Both $\mathfrak{M}_1$ and $\mathfrak{M}_2$ see their minimum expansion over $\mathcal{C}'$ on the unstable eigenvector of $\mathfrak{M}_2$. As does $\mathfrak{M}_3$ for $\varepsilon<\varepsilon^\star \approx 0.07735$, then on its own unstable eigenvector for $\varepsilon>\varepsilon^\star$. Since $\mathcal{C}'$ is expanding, each of the $\mathcal{K}_j(\varepsilon)$ are strictly greater than 1 across our parameter range.

First suppose $\Gamma_{p-1}$ lies entirely within one of the $A_j$. Each of its constituent line segments $L(x_i,v_i)$ can be defined by an end point $x_i$ and the vector $v_i$ taking $x_i$ to the other end point, with $v_i \in \mathcal{C}'$. Satisfying (C3') we let $\Gamma_p = H^{-j}(\Gamma_{p-1})$, then each $L(x_i,v_i)$ is mapped to a new segment $L(H^{-j}(x_i), \mathfrak{M}_jv_i)$ which lies in $A$, is aligned in $\mathcal{C}'$ and has expanded in diameter by a factor of at least $\mathcal{K}_j(\varepsilon)$.

As the union of these new line segments, $\Gamma_p$ satisfies (C1') and (C2'). It does not double back on itself since $\mathfrak{M}_j$ will have the same orientation preserving (or reversing) effect on each of the $v_i$. This satisfies (C0') and tells us that the diameter of $\Gamma_p$ is the sum of the diameters of the new line segments\footnote{Assuming its not 1, at which point $\Gamma_p$ has non-simple intersection with some $A_j$}, meaning its diameter has expanded by at least the factor $\mathcal{K}_j(\varepsilon)$, satisfying (C4').

The above is the simplest case we will consider. The picture becomes more complicated as we allow intersections with multiple $A_j$. First assume that $\Gamma_{p-1}$ intersects $A_1$ and one of $A_2$ or $A_3$. We proceed by restricting to one of the $A_j$, $\Gamma^j := \Gamma_{p-1} \cap A_j$, and expanding from there, $\Gamma_p = H^{-j}(\Gamma^j)$. By the same reasoning given for the $\eta$-perturbation, we require
\begin{equation}
    \label{eq:K2ineq}
    \mathcal{K}_2(\varepsilon)  > \frac{1}{1-\frac{1}{\mathcal{K}_1(\varepsilon)}}
\end{equation}
and
\begin{equation}
    \label{eq:K3ineq}
    \mathcal{K}_3(\varepsilon)  > \frac{1}{1-\frac{1}{\mathcal{K}_1(\varepsilon)}}.
\end{equation}
Solving (\ref{eq:K2ineq}) gives $\varepsilon>\varepsilon_0 \approx 0.00925$, the lower bound on our parameter range. Solving (\ref{eq:K3ineq}) gives $\varepsilon< \varepsilon_3 \approx 0.0885$, slightly larger than the upper bound on our parameter range $\varepsilon_2$.

Next assume that $\Gamma_p$ intersects $A_1$, $A_2$, and $A_3$. The case where $\Gamma_p$ intersects $A_2$ and $A_3$ but not $A_1$ follows as a trivial consequence and will be addressed at the end of the proof. Clearly if the proportion of the diameter in $A_1$ exceeds $\mathcal{K}_1(\varepsilon)^{-1}$,
\[\frac{\mathrm{diam}(\Gamma^1)}{\mathrm{diam}(\Gamma_{p-1})} >  \frac{1}{ \mathcal{K}_1(\varepsilon)},\]
then we can take $\Gamma_p = H^{-1}(\Gamma^1)$ to satisfy (C0'-5'). Otherwise we have to expand from some subset of $\Gamma^2 \cup \Gamma^3$, giving $\Gamma_p$ such that
\[ \mathrm{diam}(\Gamma_p) > \frac{1}{1-\frac{1}{\mathcal{K}_1(\varepsilon)}} \, \mathrm{diam}(\Gamma^2 \cup \Gamma^3).\]
To reduce the $\varepsilon$ dependence of the problem and simplify the equations, we will take take
\[ c = \sup_{\varepsilon_0 < \varepsilon \leq \varepsilon_2} \frac{1}{1-\frac{1}{\mathcal{K}_1(\varepsilon)}} \approx 1.4765 \]
and show
\begin{equation}
\label{eq:expansCondition}
    \mathrm{diam}(\Gamma_p) > c \, \mathrm{diam}(\Gamma^2 \cup \Gamma^3).
\end{equation}
We will give an argument for expanding $\Gamma_{p-1}$ which intersects the lower portion of $A_2$. The argument for the upper portion is entirely analogous due to the $180^\circ$ rotational symmetry of both the partition of $A$ and the invariant cone.

Consider the subdivision of $A_2$ into points which remain in $A$ for a further iteration of $H^{-1}$ after returning, $A_4$, and those which do not, $A_5$. This subdivision is shown in Figure \ref{fig:greenRed}. The points labelled points are
\[ Q_1 = \left( \frac{- 4 \varepsilon^{3} - 2 \varepsilon^{2} + \varepsilon + \frac{1}{2}}{12 \varepsilon^{2} + 16 \varepsilon + 1} , 0 \right) \quad \text{and} \quad Q_2 = \left( \frac{1+2\varepsilon}{2+2\varepsilon}, \frac{3\varepsilon+2\varepsilon^2}{2+2\varepsilon} \right) \]
so that the segment $L_1$ along the $A_4$, $A_5$ boundary has gradient
\[ k_1 = \frac{12\varepsilon^2 + 16\varepsilon+1}{(2\varepsilon+1)(2\varepsilon+5)}. \]
The segment along the $A_4$, $A_3$ boundary has gradient
\[ k_2 = \frac{4\varepsilon}{2\varepsilon+1}. \]
Strictly speaking, at larger $\varepsilon$ values $A_4$ contains an additional region in the lower part of $A_5$ near $\left(\frac{1}{2}-\varepsilon,0\right)$. The only assumption we make about points in $A_5$ is that they return to $A$ after two iterations, so treating this additional region as part of $A_5$ has no impact on our analysis.

The region $A_4$ has some useful properties. Firstly, like $A_3$, segments contained within $A_4$ return to $A$ after 3 iterations. This\footnote{Together with the fact that $\mathfrak{M}_3$ and $\mathfrak{M}_4$ have the same orientation reversing effect on the invariant cone} means we can take $\Gamma_p = H^{-3}(\Gamma_3 \cup \Gamma_4)$ and have a much larger initial curve to expand from. Secondly, diameter expansion is generally strong from $A_4$. The itinerary path is $\mathfrak{b}\mathfrak{a}\mathfrak{a}$ with corresponding matrix 
\[\mathfrak{M}_4 = DH_1^{-1} DH_1^{-1} DH_0^{-1} \]
which expands vectors at least as much as any of the other $\mathfrak{M}_j$: $\mathcal{K}_4(\varepsilon) > \mathcal{K}_j(\varepsilon)$ for all $\varepsilon_0 < \varepsilon \leq \varepsilon_2$,  $j=1,2,3$.
Finally if $\Gamma_{p-1}$ intersects $A_5$, then it must traverse $A_4$ since, by assumption, it also intersects $A_3$. The case where $\Gamma_{p-1}$ does not intersect $A_5$ is trivial, reducing to the case where $\Gamma_{p-1}$ only intersects $A_1$ and $A_3$, since $A_3$ and $A_4$ both map into $A$ under $H^{-3}$ and $\mathcal{K}_4 > \mathcal{K}_3$.

Assume, then, that $\Gamma_{p-1}$ intersects $A_5$. Let $\Gamma_p = H^{-3}(\Gamma^3 \cup \Gamma^4)$. Our aim is to minimise $\mathrm{diam}(\Gamma_p)$, considering all possible curves $\Gamma_{p-1}$ dictated by the invariant cone, and showing that it still satisfies (\ref{eq:expansCondition}). To arrive at the minimal case we can make several assumptions. Firstly, $\mathrm{diam}(\Gamma^3) = 0$. The condition that we intersect $A_3$ does not stipulate any minimum diameter in $A_3$, it can be arbitrarily small. Since $\mathfrak{M}_3$ and $\mathfrak{M}_4$ have the same orientation reversing effect on vectors in the cone, assuming $\Gamma_p$ does not have diameter 1 (at which point we has non-simple intersection with some $A_j$),
\[ \mathrm{diam}(\Gamma_p) \geq \mathcal{K}_3(\varepsilon)\, \mathrm{diam}(\Gamma^3) + \mathcal{K}_4(\varepsilon)\, \mathrm{diam}(\Gamma^4). \]
Comparing with (\ref{eq:expansCondition}), taking $\mathrm{diam}(\Gamma^3)>0$ grows the RHS of (\ref{eq:expansCondition}) by $c \,\mathrm{diam}(\Gamma^3)$, but grows the LHS of (\ref{eq:expansCondition}) by at least $\mathcal{K}_3(\varepsilon)\mathrm{diam}(\Gamma^3)$. Since $\mathcal{K}_3(\varepsilon) > c$ for every $\varepsilon_0 < \varepsilon \leq \varepsilon_2$, in the minimal case $\mathrm{diam}(\Gamma^3) = 0$. We note that the condition (\ref{eq:expansCondition}) now looks like
\[    \mathrm{diam}(H^{-3}(\Gamma^4)) > c \, \mathrm{diam}(\Gamma^4 \cup \Gamma^5), \]
which is satisfied if
\begin{equation}
    \label{eq:refinedExpansCondition}
    \mathcal{K}_4(\varepsilon) > c \, \frac{ \mathrm{diam}(\Gamma^4) + \mathrm{diam}(\Gamma^5) }{ \mathrm{diam}(\Gamma^4)}.
\end{equation}
To show that this holds, we will put lower bounds on
\begin{equation}
    \label{eq:diameterProportion}
    \frac{\mathrm{diam}(\Gamma^4)}{\mathrm{diam}(\Gamma^4) + \mathrm{diam}(\Gamma^5)}
\end{equation}
and $\mathcal{K}_4(\varepsilon)$, then compare their product with $c$.

By a purely geometric argument, comparing the admissible gradients given by the invariant cone with the lines which make up the partition boundaries, we have a lower bound
\[ \frac{\mathrm{diam}(\Gamma^4)}{\mathrm{diam}(\Gamma^4) + \mathrm{diam}(\Gamma^5)} > \frac{(2\varepsilon+1)(2\varepsilon+1-2k_5^+)}{(2\varepsilon+1)(-k_4^-(2\varepsilon+3)-k_5^+(2\varepsilon+5))+12\varepsilon^2+16\varepsilon+1} := \mathcal{B}_1(\varepsilon) \]
where $k_5^+ = \sup_{\varepsilon}\mathfrak{g}_2^u(\varepsilon) \approx -0.08750$ and $k_4^- = \inf_{\varepsilon} \mathfrak{g}_3^u(\varepsilon) \approx -0.6688$. The calculation of this bound can be found in the appendix.

We will now put a lower bound on $\mathcal{K}_4(\varepsilon)$, the minimum expansion of $\mathfrak{M}_4$ over the minimal cone. This is on the anti-clockwise boundary, $v_u(\mathfrak{M}_2)$, which can be described as the vector $(1,k_5(\varepsilon))^T$ with
\[ k_5(\varepsilon) = \frac{\varepsilon - \sqrt{\varepsilon \left(4 \varepsilon^{2} + 5 \varepsilon + 1\right)}}{2 \varepsilon + 1} < 0. \] By calculating the matrix entries of $\mathfrak{M}_4$ and noting that $\mathfrak{M}_4$ reverses the orientation of vectors, one can show that
\[\mathcal{K}_4(\varepsilon) = \frac{3+46\varepsilon+52\varepsilon^2+8\varepsilon^3}{1+2\varepsilon-4\varepsilon^2-8\varepsilon^3} - \frac{12 \varepsilon + 14}{1-4\varepsilon^2} k_5(\varepsilon). \]
Let $L(\varepsilon)$ be the linear approximation for $k_5(\varepsilon)$,
\begin{equation*}
    \begin{split}
        L(\varepsilon) &= \frac{\varepsilon-\varepsilon_0}{\varepsilon_2-\varepsilon_0} (k_5(\varepsilon_2)-k_5(\varepsilon_0)) + k_5(\varepsilon_0) \\
        & = \frac{\varepsilon-\varepsilon_0}{\varepsilon_2-\varepsilon_0} (k_5^--k_5^+) + k_5^+.
    \end{split}
\end{equation*}
One can verify that $\frac{\d }{\d \varepsilon}k_5<0$ and $\frac{\d^2 }{\d \varepsilon^2}k_5>0$ for $\varepsilon_0<\varepsilon \leq \varepsilon_2$, so that $L(\varepsilon) > k_5(\varepsilon)$ across this parameter range and is equal at its extremes. This implies
\[ \mathcal{K}_4(\varepsilon) \geq \frac{3+46\varepsilon+52\varepsilon^2+8\varepsilon^3}{1+2\varepsilon-4\varepsilon^2-8\varepsilon^3} - \frac{12 \varepsilon + 14}{1-4\varepsilon^2} L(\varepsilon) := \mathcal{B}_2(\varepsilon) \]
To show condition (\ref{eq:refinedExpansCondition}), and complete this final case, it is sufficient to show that
\begin{equation}
\label{eq:finalCondition}
    \mathcal{B}_1(\varepsilon) \mathcal{B}_2(\varepsilon) > c \approx 1.4765.
\end{equation}
One can show that $\mathcal{B}_1(\varepsilon) \mathcal{B}_2(\varepsilon)$ is monotone increasing (appendix) over $\varepsilon_0<\varepsilon \leq \varepsilon_2$ and therefore takes its minimal value at $\varepsilon_0$. Plugging in this value gives
\[ \mathcal{B}_1(\varepsilon_0) \mathcal{B}_2(\varepsilon_0) \approx 1.532235, \]
which establishes (\ref{eq:finalCondition}).

The case where $\mathrm{diam}(\Gamma^1) = 0$ follows as a trivial consequence. $\mathcal{B}_1(\varepsilon)$ is still a lower bound for the proportion of $\Gamma_{p-1}$ in $A_3 \cup A_4$ so we only need to compare $\mathcal{B}_1(\varepsilon) \mathcal{B}_2(\varepsilon)$ against $c=1$ in this case.
\end{proof}

One can follow an entirely analogous argument to prove the equivalent lemma for growth in forwards time:

\begin{lemma}
\label{lemma:epsUnstab}
Let $\Gamma_{p-1}$ be a piecewise linear curve satisfying
\begin{enumerate}[label={(C\arabic*)}]
    \setcounter{enumi}{-1}
    \item $\Gamma_{p-1}$ does not double back on itself,
    \item $\Gamma_{p-1} \subset \mathfrak{a}$,
    \item Each line segment in $\Gamma_{p-1}$ is aligned with some vector in the minimal invariant cone $\mathcal{C}$ for $H$,
\end{enumerate}
and which has simple intersection with each of the $\mathfrak{a}_i$. There there exists a piecewise linear curve $\Gamma_p$ satisfying (C0), (C1), (C2),
\begin{enumerate} [label={(C\arabic*)}]
    \setcounter{enumi}{2}
    \item $\Gamma_p \subset H^{i}(\Gamma_{p-1})$ for a chosen $i \in \{ 1,2,3\}$,
    \item There exists $\delta>0$ such that $\mathrm{diam}(\Gamma_p) \geq (1+\delta)\,\mathrm{diam}(\Gamma_{p-1})$,
\end{enumerate}
where we measure the diameter of a curve using its projection to the $y$-axis.
\end{lemma}

We now give the argument for mapping into $h$-segments and $v$-segments, whose definitions we generalise to piecewise linear curves which connect the relevant boundaries of $A_1$.

\begin{lemma}
\label{lemma:epsHseg}
Let $\Gamma_P \subset A$ be a piecewise linear curve with each of its line segments aligned with a vector in $\mathcal{C}'$. If $\Gamma_P$ has non-simple intersection with some $A_i$, then $H^{-k}(\Gamma_P)$ contains a $h$-segment for some $k \in \{0,4\}$.
\end{lemma}

\begin{proof}

In comparison with Lemma \ref{lemma:etaHseg}, we have fewer non-trivial cases to consider. We claim that any $\Gamma_P$ which has non-simple intersection with $A_2$ contains a $h$-segment, that is, it can only connect $A_2$ to itself by traversing $A_1$. Since if $\Gamma_P$ were to connect the two parts of $A_2$ through $A_3$, it would have to contain a segment with gradient
\[ \mathfrak{g} < \frac{\frac{1}{2}-\varepsilon-2\varepsilon}{\frac{1}{2}-\varepsilon - \left(\frac{1}{2} + \varepsilon \right)} = -\frac{1-6\varepsilon}{4\varepsilon} =: h(\varepsilon), \]
the gradient of the line segment joining the points $\left( \frac{1}{2}-\varepsilon,\frac{1}{2}-\varepsilon \right)$ and $\left( \frac{1}{2}+\varepsilon,2\varepsilon \right)$. However $\mathfrak{g}$ is bounded from below by $\mathfrak{g}_3^u(\varepsilon)$ with
\[ \mathfrak{g}_3^u(\varepsilon) \geq \mathfrak{g}_3^u(\varepsilon_2) \approx -0.6688 \]
across $\varepsilon_0<\varepsilon \leq \varepsilon_2$. Now
\[ h(\varepsilon) \leq h(\varepsilon_2) \approx -1.4397 \]
across the range, so that $\mathfrak{g}> h(\varepsilon)$ at each $\varepsilon$. Hence if $\Gamma_P$  has non-simple intersection with $A_2$, it follows that it contains a $h$-segment. The same clearly holds if $\Gamma_P$ has non-simple intersection with $A_3$. 

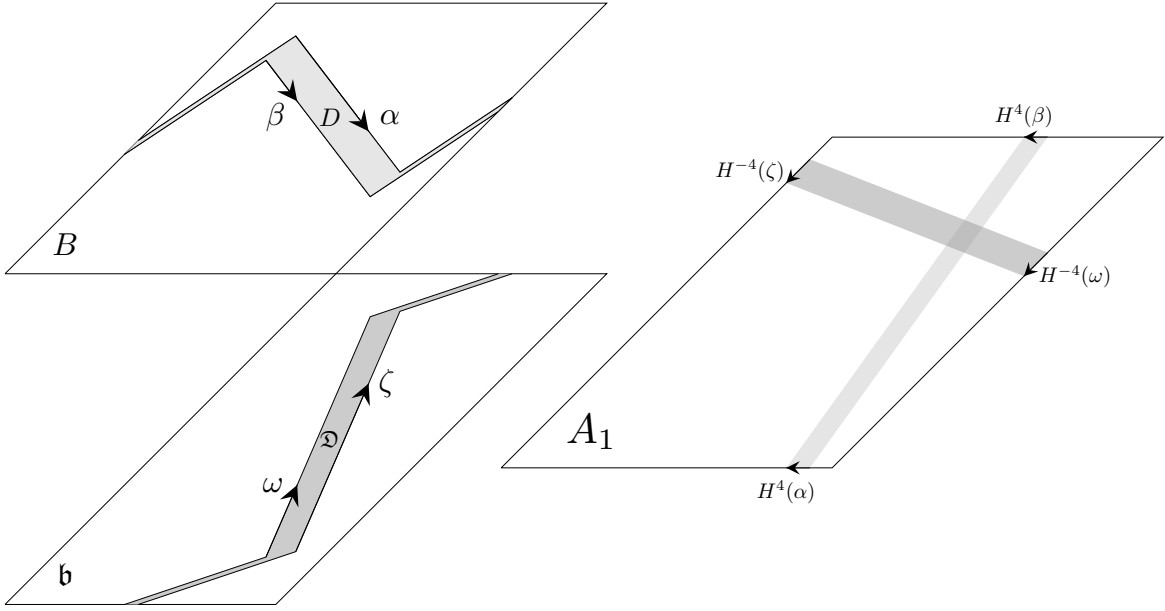
\begin{figure}
    \centering
    \begin{tikzpicture}
       
    \node[scale=0.8] at (-3,0) {
    
    \begin{tikzpicture}
       
    \draw (0,0) -- (4.5,0) -- (10,5.5) -- (5.5,5.5) -- (0,0);

    \fill[gray!40]  (1.98678163887024, 0) -- (4.33691756272402, 0.788530465949822) -- (6.06813417190776, 4.78511530398323) -- (8.19877626891103, 5.5) -- (8.42556926731037, 5.5) -- (6.56319740312437, 4.87512680056807) -- (4.83198079394063, 0.878541962534660) -- (2.21357463726958, 0);
    
    \draw (1.98678163887024, 0) -- (4.33691756272402, 0.788530465949822) -- (6.06813417190776, 4.78511530398323) -- (8.19877626891103, 5.5);
    
    \draw [decoration={markings,mark=at position 1 with
    {\arrow[scale=3,>=stealth]{>}}},postaction={decorate}] (4.33691756272402, 0.788530465949822) -- (0.7*4.33691756272402+0.3*6.06813417190776,0.7*0.788530465949822+0.3*4.78511530398323);
    
    \draw [decoration={markings,mark=at position 1 with
    {\arrow[scale=3,>=stealth]{>}}},postaction={decorate}] (4.83198079394063, 0.878541962534660) -- (0.3*4.83198079394063+0.7*6.56319740312437,0.3*0.878541962534660+0.7*4.87512680056807);
    
    \draw (8.42556926731037, 5.5) -- (6.56319740312437, 4.87512680056807) -- (4.83198079394063, 0.878541962534660) -- (2.21357463726958, 0);
    
    \node[scale=1.5] at (0.7*4.33691756272402+0.3*6.06813417190776-0.4,0.7*0.788530465949822+0.3*4.78511530398323) {$\omega$};
    
    \node[scale=1.5] at (0.3*4.83198079394063+0.7*6.56319740312437+0.3,0.3*0.878541962534660+0.7*4.87512680056807) {$\zeta$};
    
    \node[scale=1.1] at (5.4,2.75) {$\mathfrak{D}$}; 

    \node[scale=1.5] at (1,0.5) {$\mathfrak{b}$};

    % alphabeta
    
    \node[scale=1.5] at (1,6) {$B$};
    
    \draw (0,5.5) -- (5.5,5.5) -- (10,10) -- (4.5,10) -- (0,5.5);
    
    \draw (2.21357463726958, 7.71357463726958)  -- (4.83198079394063, 9.45343883140597) -- (6.56319740312437, 7.18807060255630) -- (8.42556926731037, 8.42556926731037);

    \draw (1.98678163887024, 7.48678163887024)  --  (4.33691756272402, 9.04838709677419) -- (6.06813417190776, 6.78301886792453)  -- (8.19877626891103, 8.19877626891103) ;
    
    \filldraw[fill=gray!20] (2.21357463726958, 7.71357463726958)  -- (4.83198079394063, 9.45343883140597) -- (6.56319740312437, 7.18807060255630) -- (8.42556926731037, 8.42556926731037) --  (8.19877626891103, 8.19877626891103) -- (6.06813417190776, 6.78301886792453) -- (4.33691756272402, 9.04838709677419) -- (1.98678163887024, 7.48678163887024);

    \draw [decoration={markings,mark=at position 1 with
    {\arrow[scale=3,>=stealth]{>}}},postaction={decorate}] (4.83198079394063, 9.45343883140597) -- (0.3*4.83198079394063 + 0.7*6.56319740312437, 0.3*9.45343883140597 + 0.7*7.18807060255630);
    
    \draw [decoration={markings,mark=at position 1 with
    {\arrow[scale=3,>=stealth]{>}}},postaction={decorate}] (4.33691756272402, 9.04838709677419) -- (0.3*6.06813417190776 + 0.7*4.33691756272402, 0.3*6.78301886792453 + 0.7*9.04838709677419);
    
    \node[scale=1.5] at (4.5,8.1) {$\beta$};
    
    \node[scale=1.5] at (6.4,8.1) {$\alpha$};
    
    \node[scale=1.2] at (5.4,8.1) {$D$};
    
    \end{tikzpicture} 
    };
    
    \node[scale=0.8] at (4,0) {
    
    \begin{tikzpicture}
     
    \draw (4.5,0) -- (10,0) -- (15.5,5.5) -- (10,5.5) -- (4.5,0) ;

    %% alpha beta
    \fill[gray!40, opacity=0.5] (9.63117502280512, 0) -- (9.23777559766613, 0) -- (10+3.18750279685374, 5.5) -- (10+3.58090222199275, 5.5);
    
    \draw [decoration={markings,mark=at position 1 with
    {\arrow[scale=2,>=stealth]{>}}},postaction={decorate}] (9.63117502280512, 0) -- (9.23777559766613, 0);
    
    \draw [decoration={markings,mark=at position 1 with
    {\arrow[scale=2,>=stealth]{>}}},postaction={decorate}] (10+3.58090222199275, 5.5) -- (10+3.18750279685374, 5.5);
    
    \node at (9.23777559766613, -0.4) {$H^4(\alpha)$};
    \node at (10+3.18750279685374, 5.85) {$H^4(\beta)$};

    %% omega zeta
    \fill[gray!80, opacity=0.5] (10+3.58090222199274, 3.58090222199273) -- (10+3.18750279685374, 3.18750279685375) -- (9.23777559766615, 4.73777559766613) -- (9.63117502280512, 5.13117502280511) -- (10+3.58090222199274, 3.58090222199273);
    
    \draw [decoration={markings,mark=at position 1 with
    {\arrow[scale=2,>=stealth]{>}}},postaction={decorate}] (10+3.58090222199274, 3.58090222199273) -- (10+3.18750279685374, 3.18750279685375);
    
    \draw [decoration={markings,mark=at position 1 with
    {\arrow[scale=2,>=stealth]{>}}},postaction={decorate}]  (9.63117502280512, 5.13117502280511) -- (9.23777559766615, 4.73777559766613);

    \node at (9.23777559766615-0.6, 4.73777559766613+0.2) {$H^{-4}(\zeta)$};
    
    \node at (10+3.18750279685374+0.85, 3.18750279685375) {$H^{-4}(\omega)$};
    
    \node[scale=2] at (6,0.6) {$A_1$};
    
    \end{tikzpicture}    
    };
    
    \end{tikzpicture}

    \caption{Left: Two regions $\mathfrak{D} \subset \mathfrak{b}$ and $D \subset B$, bounded by the piecewise linear curves $\omega,\zeta$ and $\alpha,\beta$ respectively. Right: Their images in $A_1$ under $H^{-4}$ and $H^4$ respectively, establishing $h$- and $v$-segments.}
    \label{fig:finalSegments}
\end{figure}

Assume, then, that $\Gamma_P$ has non-simple intersection with $A_1$. This implies that $\Gamma_P$ connects the two sloping boundaries of $A_1$ through $\mathfrak{b} = A_2 \cup A_3$. We will show that $H^{-4}(\Gamma_P)$ contains a $h$-segment. Figure \ref{fig:finalSegments} shows a region $\mathfrak{D}\subset \mathfrak{b}$, bounded by two piecewise linear curves $\omega,\zeta$. These curves can be defined by their end points on $\partial\mathfrak{b}$ and their turning points, whose full coordinates will be given in as supplementary material. Label these points as $\omega_j,\zeta_j$, $j=1,2,3,4$ so that the $x$-coordinate increases with $j$. These curves (and hence $\mathfrak{D}$) are contained within $\mathfrak{b}$ for $\varepsilon \leq \varepsilon_2$, with $\zeta_2$ limiting onto the right boundary of $\mathfrak{b}$ ($y=x-\frac{1}{2}+\varepsilon$) as $\varepsilon\rightarrow \varepsilon_2$. In particular $\varepsilon_2 \approx 0.08504$ is the positive solution to the cubic equation
\[ 8\varepsilon^3+20\varepsilon^2+10\varepsilon-1 = 0. \]

The argument for mapping into $h$-segments is roughly analogous to that given for the $\eta$-perturbation. Applying $H^{-4}$ to $\mathfrak{D}$ gives a quadrilateral in $A_1$ with sides on its left and right boundaries (the images of $\zeta$ and $\omega$ under $H^{-4}$). Clearly any $\Gamma_P$ which joins the left and right sides of $\mathfrak{b}$ must join $\omega$ and $\zeta$ through $\mathfrak{D}$. Let $\Gamma$ be this part of the curve, then $H^{-4}(\Gamma)$ must be a piecewise linear curve joining $H^{-4}(\omega)$ and $H^{-4}(\zeta)$ through $H^{-4}(\mathfrak{D})$. That is, $H^{-4}(\Gamma)$ is a $h$-segment.

\end{proof}

\begin{lemma}
\label{lemma:epsVseg}
Let $\Gamma_P \subset \mathfrak{a}$ be a piecewise linear curve with each of its line segments aligned with a vector in $\mathcal{C}$. If $\Gamma_P$ has non-simple intersection with some $\mathfrak{a}_i$, then $H^{k}(\Gamma_P)$ contains a $v$-segment for some $k \in \{0,4\}$.
\end{lemma}

\begin{proof}
Analogous to the previous lemma, non-simple intersection with $\mathfrak{a}_2$ or $\mathfrak{a}_3$ imply that $\Gamma_P$ already contains a $v$-segment. To see this, note that if $\Gamma_P$ connected the two parts of $\mathfrak{a}_2$ through $\mathfrak{a}_3$, it would have to contain a segment with gradient
\[ g(\varepsilon) > \frac{\frac{1}{2}-\varepsilon}{2\varepsilon} =: h(\varepsilon). \]
However $g(\varepsilon)$ is bounded from above by the anti-clockwise invariant cone boundary $g_3^u(\varepsilon)$ and
\[ g_3^u(\varepsilon) \leq g_3^u(\varepsilon_2) \approx 1.669 < 2.440 \approx h(\varepsilon_2) \leq h(\varepsilon) \]
across $\varepsilon_0<\varepsilon\leq \varepsilon_2$. As before, then, it remains to assess the case where $\Gamma_P$ has non-simple intersection with $\mathfrak{a}_1$. It follows that $\Gamma_P$ joins the upper and lower boundaries of $B$ through $B$. Figure \ref{fig:finalSegments} shows a region $D$ bounded by $\partial B$ and two piecewise linear curves $\alpha$, $\beta$. These curves are contained within $B$ across $\varepsilon_0<\varepsilon\leq \varepsilon_2$, with $\alpha_2$ limiting onto the line $y=1$ as $\varepsilon \rightarrow \varepsilon_2$. Applying $H^4$ to $D$ gives a quadrilateral spanning across $A_1$, with sides $H^4(\alpha)$, $H^4(\beta)$ on its lower and upper boundaries respectively. Clearly $\Gamma_P$ must connect $\beta$ to $\alpha$ through $D$, and therefore $H^4(\Gamma_P)$ contains a $v$-segment.
\end{proof}

We are now ready to establish ergodicity over $\varepsilon_0<\varepsilon<\varepsilon_2$.

\begin{proof}[Proof of Proposition \ref{prop:epsM}]
By the same argument given in the proof of Proposition \ref{prop:etaM}, by Lemmas \ref{lemma:alignment}, \ref{lemma:epsUnstab}, \ref{lemma:epsVseg}, given any $z \in X'$ we can find $m$ such that $H^m(\gamma_u(z))$ contains a $v$-segment. Similarly by Lemmas \ref{lemma:alignment}, \ref{lemma:epsStab}, \ref{lemma:epsHseg}, given any $z' \in X'$ we can find $n$ such that $H^{-n}(\gamma_s(z'))$ contains a $h$-segment. It follows that they intersect which, since $z$ and $z'$ were arbitrary, establishes \textbf{(M)}.
\end{proof}

%% file: sections/epsPert/epsRepMan.tex
\begin{prop}
\label{prop:epsMR}
Condition \textbf{(MR)} holds for $H$ when $\varepsilon_0 < \varepsilon \leq \varepsilon_2$.
\end{prop}

\begin{proof}
Follow the same argument given in the proof of Proposition \ref{prop:etaMR}, replacing $\eta$ by $\frac{1}{2}-\varepsilon$.
\end{proof}

%% file: sections/thmProof.tex
We are now ready to prove the main theorem.
\begin{proof}[Proof of Theorem \ref{thm:Bernoulli}]

Noting that \textbf{(KS1)} and \textbf{(KS2)} were trivially satisfied, the Bernoulli property holds for $H$ over $0<\eta<\eta_1$ by Theorem \ref{thm:katok-strelcyn} and Propositions \ref{prop:etaLyp}, \ref{prop:etaM}, and \ref{prop:etaMR}. Let $\eta_2 = \frac{1}{2}-\varepsilon_2$ and $\eta_3 = \frac{1}{2}-\varepsilon_0$. Then $H$ is also Bernoulli over $\eta_2 \leq \eta<\eta_3$ by Theorem \ref{thm:katok-strelcyn} and Propositions \ref{prop:epsLyp}, \ref{prop:epsM}, and \ref{prop:epsMR}.
\end{proof}

%% file: sections/finalRemarks.tex
In summary, over the parameter range $0<\eta<\frac{1}{2}$ we have given two windows within which we can prove global hyperbolicity and two subsets where mixing results can be established. A natural question is whether these are the largest sets in which these properties hold. For hyperbolicity the bounds appear optimal, with island structures developing around period 3 orbits when $\frac{1}{3}<\eta<\frac{1}{2}-\varepsilon_1$. The itinerary for these orbits (and some neighbourhood around them) is $BCA\,BCA\,BCA\dots$ so stretching behaviour is determined by the matrix $M_3$, which is non-hyperbolic. For the mixing property, the parameter limits given are not optimal. For example, $\varepsilon_2$ is not the highest upper bound on the $\varepsilon$-mixing window that our analysis allows for, but it is very close. By considering a 5-iterate mapping into $h$- and $v$-segments, this bound could be increased only very slightly. Improving the bound $\mathcal{B}_1(\varepsilon)$ would increase it further, but would in turn complicate the already lengthy algebraic manipulations.

When following the \citeauthor{katok_invariant_1986} approach, it is typical to be left with parameter ranges where non-uniform hyperbolicity can be established but proving the mixing property is more challenging. See for example the families of maps studied in \cite{przytycki_ergodicity_1983,wojtkowski1981model}. In both of these examples, the strength of the shears is increased to break up elliptic islands and ensure an invariant cone. Indeed, the \cite{wojtkowski1981model,bullett1986invariant} map at parameter value $K=4$ exhibits similar dynamics to a variation of \citeauthor{cerbelli_continuous_2005}'s map with a double strength non-monotonic shear, i.e taking $H=G\circ F^2$. In contrast, for the perturbation considered in this work the shear strength is not varied, in particular $\int_0^1 f(y) \d y$ is independent of $\eta$.

The cornerstone of our method was establishing a partition of returns and constructing an invariant expanding cone, both to prove non-zero Lyapunov exponents and as a basis for understanding how images of local manifolds grow in diameter. This approach seems viable for proving mixing properties in other systems. For example, consider the variation of \citeauthor{cerbelli_continuous_2005}'s map, perturbing the second shear by $G(x,y) = (x,y+(1+\delta)x) \text{ mod 1}$. Non-zero Lyapunov exponents can be established by our method for $0<\delta<\delta_1 \approx 0.281$, but proving \textbf{(M)} is more challenging, largely due to the map's discontinuity cutting up the images of local manifolds.

Towards the goal of more closely resembling realistic fluid velocity profiles, natural extensions to this work include introducing non-monotonicity to the second shear and studying smooth perturbations. Both of these increase the number of derivative matrices acting on the system, which complicates the analysis. The first of these is addressed in \cite{myers_hill_exponential_2021}, taking $G$ similar to $F$ in the present article. The second is considerably more challenging and is the subject of ongoing work.

%% file: sections/appendix/epsBoundsAnalysis.tex
\subsection{Establishing the lower bound $\mathcal{B}_1(\varepsilon)$}
\begin{figure}
    \centering
     \subfigure[][]{
      \label{fig:minimisationLines}
    \begin{tikzpicture}[scale = 2.5]
   % \draw[fill=gray!20] (0,0) -- (5.5,5.5) -- (10,5.5) -- (4.5,0) -- (0,0);

   \draw (5.5,1) -- (4.5,0) -- (0,0);

    \draw (0,0) -- (363/122,0) -- (110/21,110/21-4.5) -- (5.5,1) ;

    \draw (6.2,1.7) -- (5.5,1) -- (0,0);
    
    %\draw[dashed] (110/21,110/21-4.5) -- (6.2,110/21-4.5);
    
    \node at (110/21+0.2,110/21-4.55) {$Q_2$};
    
    \draw (4.5,0) -- (6.2,0);
    
    \node at (363/122,-0.2) {$Q_1$};
     %\node at (3.1,1.2) {$A_3$};
    \node at (3.1,0.75) {$L_2$};
    \node at (6,1.7) {$L_3$};
     \node at (3.3,0.25) {$L_1$};
    % \node at (8.5,2) {$A_1$};
     %\node[white] at (3.5,0.4) {$A_4$};
    % \node at (4.25,0.2) {$A_5$};
    % \node at (6,5.1) {$A_2$};
    \end{tikzpicture} }
    \subfigure[][]{
    \label{fig:boundingSegments}
    
     \begin{tikzpicture}[scale = 6]
   % \draw[fill=gray!20] (0,0) -- (5.5,5.5) -- (10,5.5) -- (4.5,0) -- (0,0);
    \begin{scope}
    \clip(2.5,0) rectangle (4.8,2);

   \draw (5.5,1) -- (4.5,0) --(363/122,0);

    \draw (363/122-0.5,0) -- (363/122,0) -- (110/21,110/21-4.5) -- (5.5,1) -- (363/122-0.5, 302/671);

    \draw (5.7,1.2) -- (5.5,1);
    
    % Cone from red corner
    \draw[ultra thick] (4.5,0) -- (3.52279782571869,0.17855971854994702) -- (2.7938247028878846,0.5079681277977972);
    \draw (3.52279782571869,0.17855971854994702) -- (2.97972625485984, 0.541768409974517);
    
    \draw[dashed] (4.5,0) -- (3.297874945791883, 0.10518915344013313);
    
    \draw[dashed] (3.297874945791883, 0.10518915344013313) -- (2.71662487574505 , 0.493931795590009);
    
    \node at (3.5524038542766024, 0.20452746048595066+0.03) {$Q_3$};
    
    \node at (3.297874945791883, 0.10518915344013313+0.05) {$Q_3'$};
    
    \node at (3.2,0.44) {$S_4$};
    \node at (3.45,0.045) {$S_5'$};
    \node at (4,0.14) {$S_5$};
    \node at (3.04,0.22) {$S_4'$};
    
    %\draw[dashed] (110/21,110/21-4.5) -- (6.2,110/21-4.5);
    
    %\node at (110/21+0.2,110/21-4.55) {$Q_2$};
    
    \draw (4.5,0) -- (5.7,0);
    
    %\node at (363/122,-0.2) {$Q_1$};
 \end{scope}
    \node at (4.85,0.9) {$L_2$};
    \node at (4.85,0.3) {$L_3$};
    \node at (4.85,0.6) {$L_1$};

    \end{tikzpicture}
    }
    \caption{A close-up on the lower portion of $A_2$, $\varepsilon=0.05$. Part (a) shows the lines which bound the regions $A_4$ and $A_5$. Part (b) shows the curve (thickest line) across $A_2$ which minimises (\ref{eq:diameterProportion}), crossing $L_1$ at $Q_3$. Also shown is the segments $S_4$ which provides a lower bound for its diameter in $A_4$. Segments $S_4'$ and $S_5'$ are defined to give further bound on (\ref{eq:diameterProportion}) with minimal $\varepsilon$ dependence.}
 
\end{figure}
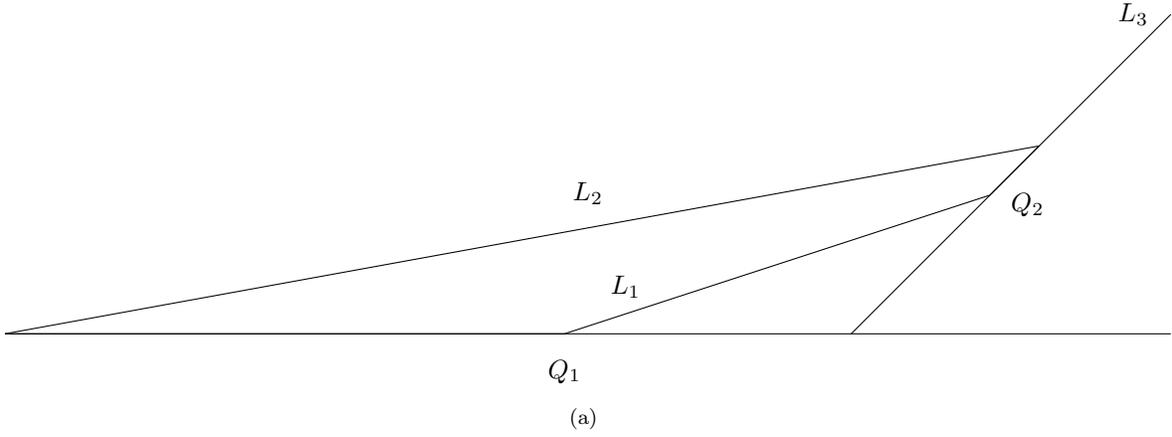
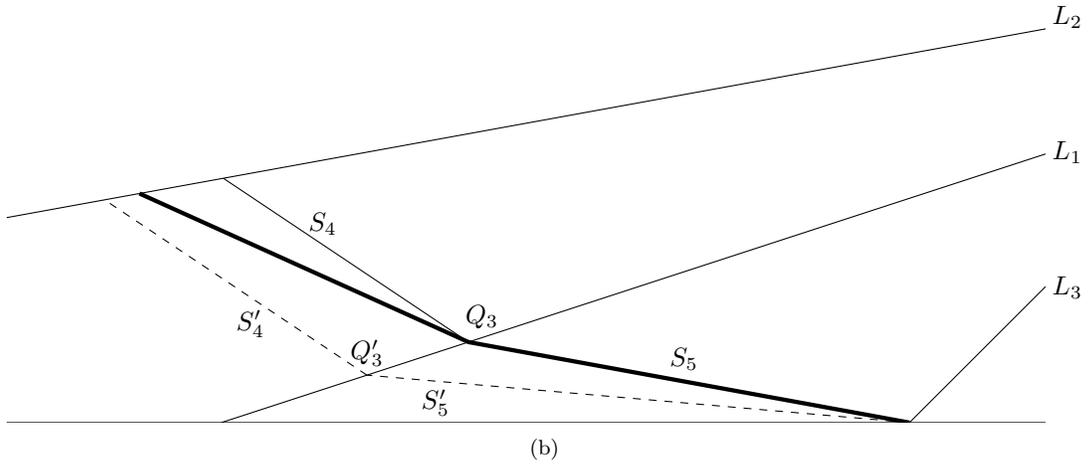
In this section we derive a lower bound
\[\frac{\mathrm{diam}(\Gamma^4)}{\mathrm{diam}(\Gamma^4) + \mathrm{diam}(\Gamma^5)} > \mathcal{B}_1(\varepsilon)\]
on the proportions of a piecewise linear curve $\Gamma_{p-1}$, constrained by the invariant cone, in the regions $A_4$ and $A_5$. We do this by maximising $\mathrm{diam}(\Gamma^5)$ and minimising $\mathrm{diam}(\Gamma^4)$, i.e we assume that $\Gamma_{p-1}$ takes the longest possible path (in diameter) across $A_5$, and the shortest possible path across $A_4$. These are straight line segments, each aligned with one of the cone boundaries. Write the gradient of segments across $A_4$ and $A_5$ as $k_4$ and $k_5$ respectively. We now have to choose where on the $L_1$ (boundary between $A_4$ and $A_5$) $\Gamma_{p-1}$ intersects so that the proportion in $A_4$ is minimal. The lines where each segment terminates are shown in Figure \ref{fig:minimisationLines}. Note that $L_2$ is the line $y=k_2x$, and $L_3$ is the line $y=x-(\frac{1}{2}-\varepsilon)$.
The diameter of the $A_4$ segment passing through $(x_1,y_1) \in L_1$ is given by
\begin{equation}
     \label{eq:A4diameter}
     \mathrm{diam}(\Gamma^4) = x_1 - \frac{y_1-k_4x_1}{k_2-k_4} 
\end{equation}
and the diameter of the $A_5$ segment passing through $(x_1,y_1) \in L_1$ is given by
\begin{equation}
    \label{eq:A5diameter}
    \mathrm{diam}(\Gamma^5) = \frac{y_1+(\frac{1}{2}-\varepsilon)-k_5x_1}{1-k_5} - x_1,
\end{equation}
valid for $(x_1,y_1) \in L_1$ above a certain threshold. This is the point $Q_3$, defined as the intersection of $L_1$ with the line $y=k_5(x-\frac{1}{2}+\varepsilon$), the lowest point on $L_1$ such that the segment in $A_5$ still intersects $L_3 \cap A$. We claim that $Q_3$ is the point where the proportion (\ref{eq:diameterProportion}) is minimal. To see this, note that as we move along the $L_1$ from $Q_2$ to $Q_3$, both diameters grow linearly. Parameterise the path as $Q_2(1-z) + Q_3z$ for $z \in [0,1]$. Now, at each $\varepsilon$, $\mathrm{diam}(\Gamma^4)(z)$ grows like $m_4 z + c_4$ for some $m_4>0$, and $c_4>0$ the diameter of the segment in $A_4$ passing through $Q_2$. Next, $\mathrm{diam}(\Gamma^5)(z)$ grows like $m_5 z$ for some $m_5>0$ since it grows from 0.
Now
\begin{equation*}
    \begin{split}
         \frac{\mathrm{diam}(\Gamma^4)}{\mathrm{diam}(\Gamma^4) + \mathrm{diam}(\Gamma^5)} (z) & = 1- \frac{\mathrm{diam}(\Gamma^5)}{\mathrm{diam}(\Gamma^4) + \mathrm{diam}(\Gamma^5)}(z) \\
         & = 1- \frac{m_5z}{m_4z+c_4+m_5z} \\
         & = 1- \frac{1}{\frac{c_4}{m_5z} + \frac{m_4}{m_5}+1}\\
    \end{split}
\end{equation*}
which is minimal at $z=1$, so (\ref{eq:diameterProportion}) is minimal at $Q_3$. We will now derive a lower bound on (\ref{eq:diameterProportion}) which has weaker $\varepsilon$ dependence.

Figure \ref{fig:boundingSegments} shows the path through $Q_3$ in bold. Its gradient in $A_5$ is given by $k_5(\varepsilon)$, aligned with the unstable eigenvector of $\mathfrak{M}_2$. Its gradient in $A_4$ is given by $k_4(\varepsilon)$, aligned with the unstable eigenvector of $\mathfrak{M}_3$. Writing the segment in $A_5$ as $S_5$, note that
\[ \frac{\mathrm{diam}(\Gamma^4)}{\mathrm{diam}(\Gamma^4) + \mathrm{diam}(\Gamma^5)} \geq \frac{\mathrm{diam}(S_4)}{\mathrm{diam}(S_4) + \mathrm{diam}(S_5)}\]
where $S_4$ is the segment in $A_4$ connecting $Q_3$ with $L_2$, with gradient aligned with the steepest possible $k_4(\varepsilon)$ over the parameter range, $k_4^- = \inf_\varepsilon k_4(\varepsilon_2) \approx -0.6688$\footnote{The minus sign in $k_4^-$ refers to it being the clockwise bound on the invariant cone}. We have equality at $\varepsilon=\varepsilon_2$. 

Now define $S_5'$ as we did $S_5$, but aligned with the least steep gradient in the parameter range, $k_5^+ = \sup_\varepsilon k_5(\varepsilon) = k_5(\varepsilon_0) \approx -0.08750$. Write its point of intersection with $L_1$ as $Q_3'$. Note that $Q_3=Q_3'$ when $\varepsilon=\varepsilon_0$. Define $S_4'$ as having the same gradient as $S_4$, but passing through $Q_3'$.

We claim that
\begin{equation}
\label{eq:segmentBound}
    \frac{\mathrm{diam}(S_4)}{\mathrm{diam}(S_4) + \mathrm{diam}(S_5)} \geq \frac{\mathrm{diam}(S_4')}{\mathrm{diam}(S_4') + \mathrm{diam}(S_5')}
\end{equation}
with equality at $\varepsilon=\varepsilon_0$. Barring this case, note that the inequality is not immediate as both $\mathrm{diam}(S_4')>\mathrm{diam}(S_4)$ and $\mathrm{diam}(S_5')>\mathrm{diam}(S_5)$. Assume the non-trivial case $\varepsilon>\varepsilon_0$ and rewrite (\ref{eq:segmentBound}) as
\[ \frac{1}{1+\frac{\mathrm{diam}(S_5)}{\mathrm{diam}(S_4)}} > \frac{1}{1+\frac{\mathrm{diam}(S_5')}{\mathrm{diam}(S_4')}}, \]
which is equivalent to
\begin{equation}
\label{eq:simplifiedSegmentBound}
 \frac{\mathrm{diam}(S_5)}{\mathrm{diam}(S_4)} < \frac{\mathrm{diam}(S_5')}{\mathrm{diam}(S_4')}.
\end{equation}
Define the diameter differences $\Delta_i = \mathrm{diam}(S_i')-\mathrm{diam}(S_i)$ and write $Q_3$ as $(x_1,y_1)$, $Q_3'$ as $(x_1',y_1')$, and $Q_1$ as $(x_0,0)$. Then $\Delta_5 = x_1-x_1'$. We can solve the line intersection equations to show that
\begin{equation}
    \label{eq:S4}
    \begin{split}   
    \mathrm{diam}(S_4)  & = x_1 - \frac{y_1-k_4^- x_1}{k_2-k_4} \\
    & = \frac{k_2x_1-y_1}{k_2-k_4^-}
    \end{split}
\end{equation}
so that
\begin{equation}
    \label{eq:delta4}
    \begin{split}   
    \Delta_4  & = x_1 - \frac{k_2x_1'-y_1'-k_2x_1+y_1}{k_2-k_4^-} \\
    & = \frac{k_2(x_1'-x_1)+k_1(x_1-x_1')}{k_2-k_4^-} \\
    & = \frac{k_1-k_2}{k_2-k_4^-} \Delta_5.
    \end{split}
\end{equation}
We can rewrite (\ref{eq:simplifiedSegmentBound}) as
\[   \frac{\mathrm{diam}(S_5') -\Delta_5}{\mathrm{diam}(S_4')-\Delta_4} < \frac{\mathrm{diam}(S_5')}{\mathrm{diam}(S_4')}, \]
which rearranges to
\[ \frac{\Delta_4}{\Delta_5} < \frac{\mathrm{diam}(S_4')}{\mathrm{diam}(S_5')}. \]
By (\ref{eq:delta4}), (\ref{eq:S4}), and $y_1'=k_1(x_1'-x_0)$ this is
\[ \frac{k_1-k_2}{k_2-k_4^-} < \frac{\frac{k_2x_1'-k_1(x_1'-x_0)}{k_2-k_4^-}}{\frac{1}{2}-\varepsilon-x_1'}, \]
which can be simplified to $ (k_1-k_2)\left(\frac{1}{2}-\varepsilon\right) < k_1x_0$.
So (\ref{eq:simplifiedSegmentBound}) holds, provided that
\[ \frac{12\varepsilon^2 + 16\varepsilon+1}{(2\varepsilon+1)(2\varepsilon+5)} - \frac{4\varepsilon}{2\varepsilon+1} < \frac{12\varepsilon^2 + 16\varepsilon+1}{(2\varepsilon+1)(2\varepsilon+5)} \cdot \frac{- 4 \varepsilon^{3} - 2 \varepsilon^{2} + \varepsilon + \frac{1}{2}}{12 \varepsilon^{2} + 16 \varepsilon + 1}, \]
which reduces to $ 1-4\varepsilon+4\varepsilon^2 < (1+2\varepsilon)^2$,
valid for all $\varepsilon>0$. This verifies the claim, giving us a lower bound
\[ \frac{\mathrm{diam}(\Gamma^4)}{\mathrm{diam}(\Gamma^4) + \mathrm{diam}(\Gamma^5)} \geq \frac{\mathrm{diam}(S_4')}{\mathrm{diam}(S_4') + \mathrm{diam}(S_5')} = \frac{k_2x_1' -y_1'}{\left(\frac{1}{2}-\varepsilon\right)(k_2-k_4^-)-y_1'+k_4^-x_1'}. \]
Noting that $y_1'$ is very small and positive\footnote{Also noting that the numerator and denominator are both positive.}, removing it from the denominator gives a new bound
\[ \frac{\mathrm{diam}(\Gamma^4)}{\mathrm{diam}(\Gamma^4) + \mathrm{diam}(\Gamma^5)} > \frac{k_2x_1' -y_1'}{\left(\frac{1}{2}-\varepsilon\right)(k_2-k_4^-)+k_4^-x_1'} :=  \mathcal{B}_1(\varepsilon)\]
which has fewer terms to consider and is still a sufficiently strong bound for our purposes. 
% Expanding the terms $k_1$, $k_2$, $x_1'$, $y_1'$ and simplifying (see appendix) yields
% \[ \mathcal{B}_1(\varepsilon) = \frac{(2\varepsilon+1)(2\varepsilon+1-2k_5^+)}{(2\varepsilon+1)(-k_4^-(2\varepsilon+3)-k_5^+(2\varepsilon+5))+12\varepsilon^2+16\varepsilon+1}. \]

\subsection{Expanding the expression for $\mathcal{B}_1(\varepsilon)$}
We will now show the expanded form of $\mathcal{B}_1(\varepsilon)$,
\begin{equation}
\label{eq:claim1}
\frac{k_2x_1' -y_1'}{\left(\frac{1}{2}-\varepsilon\right)(k_2-k_4^-)+k_4^-x_1'} = \frac{(2\varepsilon+1)(2\varepsilon+1-2k_5^+)}{(2\varepsilon+1)(-k_4^-(2\varepsilon+3)-k_5^+(2\varepsilon+5))+12\varepsilon^2+16\varepsilon+1}.
\end{equation} 
To simplify the notation, let $x=x_1'$, $k_4 = k_4^-$ and $k_5 = k_5^+$. Then $y_1' = k_1(x-x_0)$ and we can write
\begin{equation}
\label{eq:B1start}
    \mathcal{B}_1(\varepsilon) =  \frac{(k_2-k_1)x +k_1x_0}{k_4x +\left(\frac{1}{2}-\varepsilon\right)(k_2-k_4)}
\end{equation}
Let $\varphi = 2(1+2\varepsilon)(5+2\varepsilon)$. Then $\varphi k_1 = 24\varepsilon^2+32\varepsilon+2$, $\varphi k_1 x_0 = -8\varepsilon^3-4\varepsilon^2+2\varepsilon+1$, and $\varepsilon k_2 = 8\varepsilon(2\varepsilon+5)$ so that multiplying (\ref{eq:B1start}) by $\varphi/\varphi$ yields
\begin{equation*}
\begin{split}
    \mathcal{B}_1(\varepsilon) &=  \frac{(8\varepsilon(2\varepsilon+5) - (24\varepsilon^2+32\varepsilon+2))x -8\varepsilon^3-4\varepsilon^2+2\varepsilon+1}{k_4\varphi x +\left(\frac{1}{2}-\varepsilon\right)(8\varepsilon(2\varepsilon+5)-k_4\varphi)} \\
    & = \frac{ (-8\varepsilon^2+8\varepsilon-2)x +(1-2\varepsilon)(2\varepsilon
    +1)^2 }{k_4\varphi x + (1-2\varepsilon)(2\varepsilon+5)(4\varepsilon-k_4(2\varepsilon+1))}\\
    & = \frac{ -2(1-2\varepsilon)^2x +(1-2\varepsilon)(2\varepsilon
    +1)^2 }{k_4\varphi x + (1-2\varepsilon)(2\varepsilon+5)(4\varepsilon-k_4(2\varepsilon+1))}\\
    & = \frac{ -2(1-2\varepsilon)x +(2\varepsilon 
    +1)^2 }{\frac{k_4\varphi x}{1-2\varepsilon} + (2\varepsilon+5)(4\varepsilon-k_4(2\varepsilon+1))}.
\end{split}
\end{equation*}
Now by
\begin{equation*}
    \begin{split}
        x & = \frac{k_5(\varepsilon-\frac{1}{2}) +k_1x_0}{k_1-k_5} \\
        & = \frac{-k_5(2\varepsilon+5)(1-2\varepsilon)+(1-2\varepsilon)(1+2\varepsilon)}{2(2\varepsilon+5)(k_1-k_5)},
    \end{split}
\end{equation*}
we have that
\begin{equation}
\label{eq:xone}
\frac{k_4\varphi x}{1-2\varepsilon} = k_4 (2\varepsilon+1)\frac{-k_5(2\varepsilon+5)+1+2\varepsilon}{k_1-k_5} 
\end{equation}
and
\begin{equation}
\label{eq:xtwo}
-2(1-2\varepsilon)x = \frac{k_5(1-2\varepsilon)^2}{k_1-k_5} - \frac{(1-2\varepsilon)^2(1+2\varepsilon)}{(2\varepsilon+5)(k_1-k_5)}
\end{equation}
so that
\begin{equation*}
\mathcal{B}_1(\varepsilon) = \frac{k_5(2\varepsilon-1)^2(2\varepsilon+5)-(1-2\varepsilon)^2(1+2\varepsilon)+(2\varepsilon+1)^2(2\varepsilon+5)(k_1-k_5)}{k_4(2\varepsilon+1)(2\varepsilon+5)(-k_5(2\varepsilon+5)+1+2\varepsilon) + (2\varepsilon+5)^2(k_1-k_5)(4\varepsilon-k_4(2\varepsilon+1))},
\end{equation*}
where we have substituted in (\ref{eq:xone}), (\ref{eq:xtwo}) and multiplied top and bottom by $(2\varepsilon+5)(k_1-k_5)$. Write its numerator and denominator as $N(\varepsilon)$ and $D(\varepsilon)$. Expanding the $k_1$ term,
\begin{equation*}
    \begin{split}
        N(\varepsilon) &= k_5(2\varepsilon+5)\left( (2\varepsilon-1)^2-(2\varepsilon+1)^2  \right) +(2\varepsilon+1)(12\varepsilon^2+16\varepsilon+1 - (1-2\varepsilon)^2)\\
        & = -8\varepsilon k_5 (2\varepsilon+5) + (2\varepsilon+1)(8\varepsilon^2+20\varepsilon) \\
        & = (2\varepsilon+5)(4\varepsilon(2\varepsilon+1)-8\varepsilon k_5 (2\varepsilon+5))
    \end{split}
\end{equation*}
and
\begin{equation*}
    \begin{split}
        D(\varepsilon) &= k_4(2\varepsilon+1)(2\varepsilon+5)\left( -k_5(2\varepsilon+5)+1+2\varepsilon \right) - k_5(2\varepsilon+5)^2(4\varepsilon-k_4(2\varepsilon+1)) \\ 
        & \quad+ \frac{12\varepsilon^2+16\varepsilon+1}{2\varepsilon+1} (2\varepsilon+5)\left(4\varepsilon-k_4(2\varepsilon+1)\right) \\
        & = (2\varepsilon+5) \bigg( k_4 \left[ (2\varepsilon+1)(-k_5(2\varepsilon+5)+1+2\varepsilon) +k_5(2\varepsilon+5)(2\varepsilon+1) -(12\varepsilon^2+16\varepsilon+1) \right] \\
        & \quad -4\varepsilon k_5(2\varepsilon+5) + \frac{4\varepsilon}{2\varepsilon+1}(12\varepsilon+16\varepsilon+1) \bigg).
    \end{split}
\end{equation*}
Noting that the $k_4k_5$ terms cancel and $(1+2\varepsilon)^2-(12\varepsilon^2+16\varepsilon+1) = -4\varepsilon(2\varepsilon+3)$,
\begin{equation*}
   \mathcal{B}_1(\varepsilon) = \frac{4\varepsilon(2\varepsilon+1)-8\varepsilon k_5}{-4\varepsilon k_4(2\varepsilon+3)-4\varepsilon k_5 (2\varepsilon+5) + \frac{4\varepsilon}{2\varepsilon+1}(12\varepsilon^2+16\varepsilon+1)}.
\end{equation*}
Multiplying top and bottom by $\frac{2\varepsilon+1}{4\varepsilon}$ establishes (\ref{eq:claim1}).

\subsection{$\mathcal{B}_1(\varepsilon)\mathcal{B}_2(\varepsilon)$ is monotone increasing}
Starting with the bound on $K_4(\varepsilon)$,
\begin{equation*}
    \begin{split}
        \mathcal{B}_2(\varepsilon) &= \frac{3+46\varepsilon+52\varepsilon^2+8\varepsilon^3}{1+2\varepsilon-4\varepsilon^2-8\varepsilon^3} - \frac{12 \varepsilon + 14}{1-4\varepsilon^2} L(\varepsilon) \\
        & = \frac{3+46\varepsilon+52\varepsilon^2+8\varepsilon^3}{(1-2\varepsilon)(1+2\varepsilon)^2} - \frac{12 \varepsilon + 14}{(1-2\varepsilon)(1+2\varepsilon)} L(\varepsilon) \\
        & = \frac{3+46\varepsilon+52\varepsilon^2+8\varepsilon^3-(1+2\varepsilon)(12\varepsilon+14)L(\varepsilon)}{(1-2\varepsilon)(1+2\varepsilon)^2}.
    \end{split}
\end{equation*}
Combining with our expanded expression for $\mathcal{B}_1(\varepsilon)$,
\begin{equation*}
        \mathcal{B}_1(\varepsilon)\mathcal{B}_2(\varepsilon)  = \frac{(2\varepsilon+1-2k_5^+)(  3+46\varepsilon+52\varepsilon^2+8\varepsilon^3-(1+2\varepsilon)(12\varepsilon+14)L(\varepsilon))}{(1-2\varepsilon)(1+2\varepsilon)^2(-k_4^-(2\varepsilon+3)-k_5^+(2\varepsilon+5))+ (1-2\varepsilon)(1+2\varepsilon)(12\varepsilon^2+16\varepsilon+1)}
\end{equation*}
where we have divided through by $(1+2\varepsilon)/(1+2\varepsilon)$. Write its numerator and denominator as $P(\varepsilon)$ and $Q(\varepsilon)$, then $\mathcal{B}_1(\varepsilon)\mathcal{B}_2(\varepsilon)$ is monotone increasing if $P'Q-PQ'>0$. Note that as a linear function, $L'(\varepsilon) =k_6 \approx -1.85175$ is constant. The factors derived from $P$ are then
\begin{equation*}
    \begin{split}
        P(\varepsilon) &= (2\varepsilon+1-2k_5^+)(  3+46\varepsilon+52\varepsilon^2+8\varepsilon^3-(24\varepsilon^2+40\varepsilon+14)L(\varepsilon)),\\
        P'(\varepsilon) & =  6+92\varepsilon+104\varepsilon^2+16\varepsilon^3-(48\varepsilon^2+80\varepsilon+28)L(\varepsilon) \\
        & \quad+ (2\varepsilon+1-2k_5^+)(  46+104\varepsilon+24\varepsilon^2-(48\varepsilon+40)L(\varepsilon) - k_6(24\varepsilon^2+40\varepsilon+14))
    \end{split}
\end{equation*}
which, since $(2\varepsilon+1-2k_5^+)>0$, $L(\varepsilon)<0$, and $k_6<0$, are both positive for $\varepsilon>0$. Hence over the parameter range $\varepsilon_0 <\varepsilon \leq \varepsilon_2$, $P$ is maximal at $\varepsilon_2$. Differentiating again, one can verify that $P''>0$, so that $P'$ is bounded below by $P'(\varepsilon_0)$. Now for the factors derived from the denominator,
\begin{equation*}
    \begin{split}
        Q(\varepsilon) &= (8\varepsilon^3+4\varepsilon^2-2\varepsilon_1)(k_4^-(2\varepsilon+3)+k_5^+(2\varepsilon+5))+ (1-4\varepsilon^2)(12\varepsilon^2+16\varepsilon+1),\\
        Q'(\varepsilon) & = (24\varepsilon^2+8\varepsilon-2)(k_4^-(2\varepsilon+3)+k_5^+(2\varepsilon+5))+(8\varepsilon^3+4\varepsilon^2-2\varepsilon_1)(2k_4^-+2k_5^+) \\
        & \quad -8\varepsilon(12\varepsilon^2+16\varepsilon+1)+(1-4\varepsilon^2)(24\varepsilon+16)
    \end{split}
\end{equation*}
which, since $8\varepsilon^3+4\varepsilon^2-2\varepsilon_1<0$ and $-8\varepsilon(12\varepsilon^2+16\varepsilon+1)+(1-4\varepsilon^2)(24\varepsilon+16) = 16+16\varepsilon+\dots >0$ over the parameter range, are also both positive. Hence $Q$ bounded below by $Q(\varepsilon_0)$. Again, one can verify that $Q''<0$ so that $Q'$ is bounded above by $Q'(\varepsilon_0)$.

Hence $P'Q-PQ' > P'(\varepsilon_0)Q(\varepsilon_0)-P(\varepsilon_2)Q'(\varepsilon_0) \approx 29.853$, positive as required.

%% file: sections/etaPert/supp.tex
\subsubsection{Eigenvector gradients}

Shown below are the gradients of the eigenvectors defined in section \ref{sec:etaLyp}, as functions in $\eta$. To simplify the expressions, let $a(\eta) = - 16 \eta^{3} + 44 \eta^{2} - 36 \eta + 9$ and $b(\eta) = - 36 \eta^{3} + 93 \eta^{2} - 54 \eta + 9$.

\[ g_1^u =   \frac{2 - 2 \eta}{\sqrt{5 - 4 \eta} - 1}    \quad     g_1^s =   \frac{2 \left(\eta - 1\right)}{\sqrt{5 - 4 \eta} + 1}   \]
\[ g_2^u =  \frac{2 \eta - \sqrt{a(\eta)} - 3}{4 \left(\eta - 1\right)}     \quad     g_2^s =   \frac{2 \eta + \sqrt{a(\eta)} - 3}{4 \left(\eta - 1\right)}
   \]
\[ g_3^u =    \frac{3 \eta^{2} + (1-\eta) \sqrt{b(\eta)} - 8 \eta + 3}{6 \eta^{2} - 12 \eta + 4}   \quad     g_3^s =   \frac{3 \eta^{2} + (\eta -1) \sqrt{b(\eta)} - 8 \eta + 3}{6 \eta^{2} - 12 \eta + 4}    \]

\[ \mathfrak{g}_1^u =   \frac{2 \left(\eta - 1\right)}{\sqrt{5 - 4 \eta} + 1}    \quad     \mathfrak{g}_1^s =   \frac{2 - 2 \eta}{\sqrt{5 - 4 \eta} - 1}   \]
\[ \mathfrak{g}_2^u =    \frac{2 \eta + \sqrt{a(\eta)} - 1}{4 \left(\eta - 1\right)}   \quad     \mathfrak{g}_2^s =  \frac{2 \eta - \sqrt{a(\eta)} - 1}{4 \left(\eta - 1\right)}    \]
\[ \mathfrak{g}_3^u =  \frac{3 \eta^{2} + (\eta -1)\sqrt{b(\eta)} - 4 \eta + 1 }{6 \eta^{2} - 12 \eta + 4}
     \quad   \mathfrak{g}_3^s =  \frac{3 \eta^{2} + (1-\eta)\sqrt{b(\eta)} - 4 \eta + 1 }{6 \eta^{2} - 12 \eta + 4}
   \]

\subsubsection{Coordinates $P_j$}

Shown below are coordinates for the corners which define the quadrilateral $\mathcal{Q}_3$ and their images under $H^5$.

\[ P_1 = \left( \frac{\eta \left(\eta^{4} - 9 \eta^{3} + 26 \eta^{2} - 30 \eta + 12\right)}{5 \eta^{4} - 35 \eta^{3} + 80 \eta^{2} - 72 \eta + 21}, 1-\eta \right) \quad H^5(P_1) = \left( \frac{\eta^{5} - 4 \eta^{4} - 3 \eta^{3} + 29 \eta^{2} - 36 \eta + 12}{5 \eta^{4} - 35 \eta^{3} + 80 \eta^{2} - 72 \eta + 21} , 0 \right)  \]

\[ P_2 = \left( \frac{9 \eta^{4} - 51 \eta^{3} + 101 \eta^{2} - 81 \eta + 21}{5 \eta^{4} - 35 \eta^{3} + 80 \eta^{2} - 72 \eta + 21}, 1 \right) \quad H^5(P_2) = \left(\frac{- \eta^{4} - 6 \eta^{3} + 30 \eta^{2} - 36 \eta + 12}{5 \eta^{4} - 35 \eta^{3} + 80 \eta^{2} - 72 \eta + 21} , 0 \right)  \]

\[ P_3 = \left(\frac{- \eta^{5} + 13 \eta^{4} - 57 \eta^{3} + 105 \eta^{2} - 82 \eta + 21}{5 \eta^{4} - 35 \eta^{3} + 80 \eta^{2} - 72 \eta + 21} , 1 \right) \quad H^5(P_3) = \left( \frac{- \eta^{5} + 8 \eta^{4} - 23 \eta^{3} + 30 \eta^{2} - 18 \eta + 4}{5 \eta^{4} - 35 \eta^{3} + 80 \eta^{2} - 72 \eta + 21} , 1-\eta \right)  \]

\[ P_4 = \left( \frac{\eta \left(- 5 \eta^{3} + 20 \eta^{2} - 26 \eta + 11\right)}{5 \eta^{4} - 35 \eta^{3} + 80 \eta^{2} - 72 \eta + 21}, 1-\eta \right) \quad H^5(P_4) = \left(\frac{5 \eta^{4} - 20 \eta^{3} + 29 \eta^{2} - 18 \eta + 4}{5 \eta^{4} - 35 \eta^{3} + 80 \eta^{2} - 72 \eta + 21} , 1-\eta \right)  \]

\subsubsection{Coordinates $p_j$}

\begin{center}

\begin{tikzpicture}

\node at (-4,0) {

\begin{tikzpicture}

\fill[gray!50] (1.56821745948772,1.56821745948772) -- (7.99660219550444,7.99660219550444-2.5) -- (9.10245687401986,9.10245687401986-2.5) -- (2.67407213800314,2.67407213800314);

\node at (1.56821745948772-0.1,1.56821745948772) {$p_1$};

\node at (7.99660219550444+0.2,7.99660219550444-2.5) {$p_2$};

\node at (9.10245687401986+0.1,9.10245687401986-2.5) {$p_3$};

\node at (2.67407213800314-0.15,2.67407213800314) {$p_4$};

\end{tikzpicture}
};

\draw [->] (-0.5,-0.4) -- (0.5,-0.4);

\node at (0,-0.2) {$H^{-5}$};

\node at (4,0) {

\begin{tikzpicture}

\fill[gray!50] (10+2.60690015682176+0.5,2.60690015682176+0.5) -- (10+3.64218504966023,3.64218504966023) --  (8.63565081024572,8.63565081024572-2.5) -- (7.60036591740722+0.5,7.60036591740722-2.5+0.5);

\node at (10+2.60690015682176+0.5,2.60690015682176+0.2) {$H^{-5}(p_1)$};

\node at (10+3.64218504966023+0.6,3.64218504966023) {$H^{-5}(p_2)$};

\node at (8.63565081024572,8.63565081024572-2.5+0.2) {$H^{-5}(p_3)$};

\node at (7.60036591740722+0.5-0.6,7.60036591740722-2.5+0.5) {$H^{-5}(p_4)$};

\end{tikzpicture}
};

\end{tikzpicture}
\end{center}

Shown below are $x$-coordinates for the corners which define the quadrilateral $Q_3$ and their images under $H^{-5}$. The $y$-coordinates can be deduced from $p_1,p_4,H^{-5}(p_1),H^{-5}(p_2)$ being on the line $y=x$, and $p_2,p_3,H^{-5}(p_3),H^{-5}(p_4)$ being on the line $y=x-\eta$.

\[ p_1 = \frac{\eta \left(- 4 \eta^{3} + 16 \eta^{2} - 21 \eta + 9\right)}{5 \eta^{4} - 35 \eta^{3} + 80 \eta^{2} - 72 \eta + 21} \quad H^{-5}(p_1) = \frac{6 \eta^{4} - 29 \eta^{3} + 50 \eta^{2} - 36 \eta + 9}{5 \eta^{4} - 35 \eta^{3} + 80 \eta^{2} - 72 \eta + 21}  \]

\[ p_2 = \frac{- \eta^{5} + 14 \eta^{4} - 61 \eta^{3} + 110 \eta^{2} - 84 \eta + 21}{5 \eta^{4} - 35 \eta^{3} + 80 \eta^{2} - 72 \eta + 21} \quad H^{-5}(p_2) = \frac{- \eta^{5} + 9 \eta^{4} - 32 \eta^{3} + 51 \eta^{2} - 36 \eta + 9}{5 \eta^{4} - 35 \eta^{3} + 80 \eta^{2} - 72 \eta + 21}  \]

\[ p_3 = \frac{10 \eta^{4} - 55 \eta^{3} + 106 \eta^{2} - 83 \eta + 21}{5 \eta^{4} - 35 \eta^{3} + 80 \eta^{2} - 72 \eta + 21} \quad H^{-5}(p_3) = \frac{- 15 \eta^{3} + 51 \eta^{2} - 54 \eta + 17}{5 \eta^{4} - 35 \eta^{3} + 80 \eta^{2} - 72 \eta + 21}  \]

\[ p_4 = \frac{\eta \left(\eta^{4} - 8 \eta^{3} + 22 \eta^{2} - 25 \eta + 10\right)}{5 \eta^{4} - 35 \eta^{3} + 80 \eta^{2} - 72 \eta + 21} \quad H^{-5}(p_4) = \frac{\eta^{5} - 3 \eta^{4} - 12 \eta^{3} + 50 \eta^{2} - 54 \eta + 17}{5 \eta^{4} - 35 \eta^{3} + 80 \eta^{2} - 72 \eta + 21}  \]

%% file: sections/epsPert/piecewiseLinearCurves.tex
Recall the piecewise linear curves $\alpha, \beta \subset B$, $\omega,\zeta \subset \mathfrak{b}$ given in Figure \ref{fig:finalSegments}. Labelling the endpoints and turning points by increasing $x$-coordinate, the coordinates of these points are as follows.

\begin{equation*}
    \begin{split}
 \alpha_1 & =    \left(
\frac{- 16  \varepsilon^{5} - 72 \varepsilon^{4} - 72 \varepsilon^{3} - 4 \varepsilon^{2} + 19 \varepsilon + 5.5}{80 \varepsilon^{4} + 400 \varepsilon^{3} + 560 \varepsilon^{2} + 252 \varepsilon + 15 }, \frac{64 \varepsilon^{5} + 368 \varepsilon^{4} + 688 \varepsilon^{3} + 528 \varepsilon^{2} + 160 \varepsilon + 13}{80 \varepsilon^{4} + 400 \varepsilon^{3} + 560 \varepsilon^{2} + 252 \varepsilon + 15 }
 \right) \\
\alpha_2 & =    \left(- \frac{192 \varepsilon^{4} + 960 \varepsilon^{3} + 1344 \varepsilon^{2} + 720 \varepsilon + 132} {576 \varepsilon^{3} + 1056 \varepsilon^{2} + 48 \varepsilon - 360} , \frac{16 \varepsilon^{3} + 52 \varepsilon^{2} + 48 \varepsilon + 13 }{12 \varepsilon^{2} + 28 \varepsilon + 15} 
 \right) \\
  \alpha_3 & =     \left(
\frac{40 \varepsilon^{3} + 76 \varepsilon^{2} + 22 \varepsilon - 11}{24 \varepsilon^{3} + 44 \varepsilon^{2} + 2 \varepsilon - 15}
,
\frac{4 \varepsilon^{2} + 16 \varepsilon + 11}{12 \varepsilon^{2} + 28 \varepsilon + 15}
 \right) \\
\alpha_4 & =  \left(
\frac{112 \varepsilon^{4} + 448 \varepsilon^{3} + 568 \varepsilon^{2} + 240 \varepsilon + 11}{80 \varepsilon^{4} + 400 \varepsilon^{3} + 560 \varepsilon^{2} + 252 \varepsilon + 15}
,
\frac{112 \varepsilon^{4} + 448 \varepsilon^{3} + 568 \varepsilon^{2} + 240 \varepsilon + 11}{80 \varepsilon^{4} + 400 \varepsilon^{3} + 560 \varepsilon^{2} + 252 \varepsilon + 15}
 \right) \\ 
  \beta_1 & =     \left(
\frac{- 48 \varepsilon^{4} - 64 \varepsilon^{3} - 8 \varepsilon^{2} + 16 \varepsilon + 5 }{80 \varepsilon^{4} + 400 \varepsilon^{3} + 560 \varepsilon^{2} + 252 \varepsilon + 15 }
,
\frac{80 \varepsilon^{5} + 392 \varepsilon^{4} + 696 \varepsilon^{3} + 524 \varepsilon^{2} + 157 \varepsilon + 12.5}{80 \varepsilon^{4} + 400 \varepsilon^{3} + 560 \varepsilon^{2} + 252 \varepsilon + 15 } 
 \right) \\ 
  \beta_2 & =      \left(
- \frac{16 \varepsilon^{2} + 16 \varepsilon + 4} {16 \varepsilon^{2} + 16 \varepsilon - 12} 
,
\frac{\left( \varepsilon + 1.5\right) \left(2 \varepsilon + 3 \right) - 2} {2 \varepsilon + 3} 
 \right) \\          
  \beta_3 & =    \left(
\frac{24 \varepsilon^{3} + 132 \varepsilon^{2} + 66 \varepsilon - 21 }{72 \varepsilon^{2} + 24 \varepsilon - 30} 
,
\frac{- 2 \varepsilon^{2} + 2 \varepsilon + 3.5}{6 \varepsilon + 5} 
 \right) \\ 
  \beta_4 & =    \left(
\frac{16 \varepsilon^{5} + 136 \varepsilon^{4} + 456 \varepsilon^{3} + 564 \varepsilon^{2} + 237 \varepsilon + 10.5}{80 \varepsilon^{4} + 400 \varepsilon^{3} + 560 \varepsilon^{2} + 252 \varepsilon + 15} 
,
\frac{16 \varepsilon^{5} + 136 \varepsilon^{4} + 456 \varepsilon^{3} + 564 \varepsilon^{2} + 237 \varepsilon + 10.5}{80 \varepsilon^{4} + 400 \varepsilon^{3} + 560 \varepsilon^{2} + 252 \varepsilon + 15} 
 \right) \\  
  \omega_1 & =   \left(
\frac{- 48 \varepsilon^{4} - 64 \varepsilon^{3} - 8 \varepsilon^{2} + 16 \varepsilon + 5}{80 \varepsilon^{4} + 400 \varepsilon^{3} + 560 \varepsilon^{2} + 252 \varepsilon + 15}
,
0
 \right) \\   
  \omega_2 & =  \left(
- \frac{4 \varepsilon^{2} + 4 \varepsilon + 1}{4 \varepsilon^{2} + 4 \varepsilon - 3}
,
- \frac{4 \varepsilon \left(2 \varepsilon + 1\right)}{4 \varepsilon^{2} + 4 \varepsilon - 3}
 \right) \\
  \omega_3 & =   \left(
\frac{48 \varepsilon^{3} + 264 \varepsilon^{2} + 132 \varepsilon - 42}{144 \varepsilon^{2} + 48 \varepsilon - 60 }
,
\frac{240 \varepsilon^{3} + 312 \varepsilon^{2} + 36 \varepsilon - 30}{144 \varepsilon^{2} + 48 \varepsilon - 60 }
 \right) \\    
  \omega_4 & =     \left(
\frac{16 \varepsilon^{5} + 136 \varepsilon^{4} + 456 \varepsilon^{3} + 564 \varepsilon^{2} + 237 \varepsilon + 10.5}{80 \varepsilon^{4} + 400 \varepsilon^{3} + 560 \varepsilon^{2} + 252 \varepsilon + 15 }
,
 \frac{1}{2} + \varepsilon
 \right) \\
  \zeta_1 & =    \left(
\frac{- 16 \varepsilon^{5} - 72 \varepsilon^{4} - 72 \varepsilon^{3} - 4 \varepsilon^{2} + 19 \varepsilon + 5.5}{80 \varepsilon^{4} + 400 \varepsilon^{3} + 560 \varepsilon^{2} + 252 \varepsilon + 15 }
,
0
 \right) \\  
  \zeta_2 & =   \left(
- \frac{192 \varepsilon^{4} + 960 \varepsilon^{3} + 1344 \varepsilon^{2} + 720 \varepsilon + 132 }{576 \varepsilon^{3} + 1056 \varepsilon^{2} + 48 \varepsilon - 360 }
,
- \frac{\varepsilon \left(384 \varepsilon^{3} + 1728 \varepsilon^{2} + 1824 \varepsilon + 528 \right)}{576 \varepsilon^{3} + 1056 \varepsilon^{2} + 48 \varepsilon - 360 }
 \right) \\ 
  \zeta_3 & =    \left(
\frac{960 \varepsilon^{3} + 1824 \varepsilon^{2} + 528 \varepsilon - 264}{576 \varepsilon^{3} + 1056 \varepsilon^{2} + 48 \varepsilon - 360 }
,
\frac{576 \varepsilon^{4} + 2112 \varepsilon^{3} + 1728 \varepsilon^{2} + 48 \varepsilon - 180}{576 \varepsilon^{3} + 1056 \varepsilon^{2} + 48 \varepsilon - 360 }
 \right) \\ 
  \zeta_4 & =   \left(
\frac{112 \varepsilon^{4} + 448 \varepsilon^{3} + 568 \varepsilon^{2} + 240 \varepsilon + 11 }{80 \varepsilon^{4} + 400 \varepsilon^{3} + 560 \varepsilon^{2} + 252 \varepsilon + 15 }
,
 \frac{1}{2} + \varepsilon
 \right) \\
    \end{split}
\end{equation*}